\documentclass{amsart}
\title{The enumerative geometry of $K3$ surfaces and modular forms}
\author{Jim Bryan\\
Naichung Conan Leung}
\date{\today}
\address{Department of Mathematics\\
University of California\\
Berkeley, CA 94720}
\address{School of Mathematics\\
University of Minnesota\\
Minneapolis, MN 55455}

\usepackage{amscd,eepic,epic}
\newtheorem{thm}{Theorem}[section]
\newtheorem{theorem}[thm]{Theorem}

\newtheorem{lem}[thm]{Lemma}

\newtheorem{prop}[thm]{Proposition}
\newtheorem{proposition}[thm]{Proposition}

\newtheorem{definition}[thm]{Definition}

\newtheorem{rem1}[thm]{Remark}
\newenvironment{rem}{\begin{rem1}\em}{\end{rem1}}

\newtheorem{exer1}[thm]{Exercise}

\newcommand{\cnums} {{\mathbf C}}          % complex numbers
\newcommand{\znums} {{\mathbf Z}}		% integers
\newcommand{\qnums} {{\mathbf Q}}		% rationals

\newcommand{\im}{\operatorname{Im}}

\newcommand{\ie}{{\em i.e. }}

\newcommand{\point}{\text{pt.}}
\newcommand{\til}[1]{\Tilde{#1}}

\renewcommand{\P}{\mathbf{P}}
\newcommand{\M}{\mathcal{M}}

\begin{document}
\begin{abstract}
Let $X$ be a $K3$ surface and $C$ be a holomorphic curve in $X$ representing
a primitive homology class.  We count the number of
curves of geometric genus $g$ with $n$ nodes passing through $g$ generic
points in $X$ in the linear system $\left| C\right| $ for any $g$ and
$n$ satisfying $C\cdot C=2g+2n-2$.

When $g=0$, this
coincides with the enumerative problem 
studied by Yau and Zaslow who obtained a conjectural generating function
for the numbers. Recently, G\"ottsche has generalized their conjecture to
arbitrary $g$ in terms of quasi-modular forms. We prove these formulas
using Gro\-mov-Wit\-ten invariants for families, a degeneration argument, and
an obstruction bundle computation. 
Our methods also
apply to $\mathbf{P}^{2}$ blown up at 9 
points where we show that the ordinary
Gro\-mov-Wit\-ten invariants of genus $g$ constrained to $g$ points are also
given in terms of quasi-modular forms. 
\end{abstract}
\thanks{The first author is supported by a grant from the Ford Foundation
and the second author is supported by a NSF grant DMS-9626689.} 
\maketitle
\markboth{Enumerative geometry of $K3$}
{Enumerative geometry of $K3$}
\renewcommand{\sectionmark}[1]{}
\tableofcontents
\pagebreak

\section{Introduction}

Let $X$ be a $K3$ surface and $C$ be a holomorphic curve in $X$ representing
a primitive homology class. For any $g$ and $n$ satisfying $C\cdot
C=2g+2n-2,$ we define an invariant $N_{g}(n)$ which counts the number of
curves of geometric genus $g$ with $n$ nodes passing through $g$ generic
points in $X$ in the linear system $\left| C\right| .$ For each $g,$
consider the generating function         
\[
F_g\left( q\right) =\sum_{n=0}^\infty N_{g}(n)q^n. \]

Our main theorem gives explicit formulas for $F_{g}$ in terms of
quasi-modular forms:

\begin{theorem}[Main Theorem]\label{thm: main thm}
For any $g,$ we have
\begin{eqnarray*}
F_g\left( q\right) &=&\left(\sum_{k=1}^{\infty
}k(\sum_{d|k}d)q^{k-1} \right)^{g} \prod _{m=1}^{\infty }(1-q^{m})^{-24}\\
&=&\left(\frac{d}{dq}G_{2}(q) \right)^{g}\frac{q}{\Delta (q)}.
\end{eqnarray*}  
\end{theorem}

So for example, we have
\begin{eqnarray*}
F_{0}&=&1+ 24q+ 324q^{2}+3200q^{3}+\cdots \\
F_{1}&=&1+ 30q+ 480q^{2}+5460q^{3}+\cdots \\
F_{2}&=&1+36q+672q^{2}+8728q^{3}+\cdots \\
F_{3}&=&1+42q+900q^{2}+13220q^{3}+\cdots .
\end{eqnarray*}

If we write $q=e^{2\pi i\tau }$ then $\Delta \left( \tau \right)
=q\prod_{m=1}^{\infty }\left( 1-q^{m}\right) ^{24}=\eta \left( \tau \right)
^{24}$ is a modular form of weight $12$ where $\eta \left( \tau \right) $
is the Dedekind $\eta $ function. In particular, for any $\left(
\begin{smallmatrix}
a & b \\
c & d
\end{smallmatrix}
\right) \in SL\left( 2,\mathbf{Z}\right) $ and $\im (\tau) >0,$ we have \[
\Delta \left( \frac{a\tau +b}{c\tau +d}\right) =\left( c\tau +d\right)
^{12}\Delta \left( \tau \right) . 
\]
$G_{2}(q)$ is the Eisenstein series
$$
G_{2}(q)=\frac{-1}{24}+\sum_{k=1}^{\infty }\sigma (k)q^{k}
$$
where $\sigma (k)=\sum_{d|k}d$. $G_{2}$ and its derivatives are
quasi-modular forms \cite{Gott-conj}.

One crucial property for us is the fact that $\prod_{m=1}^{\infty
}\left( 1-q^{m}\right) ^{-1}$ is the generating function of the partition
function $p\left( d\right) $. Namely         
\begin{eqnarray*}
\prod_{m=1}^{\infty }\left( 1-q^{m}\right) ^{-1} &=&\sum_{d=0}^{\infty
}p\left( d\right) q^{d} \\ 
&=&\allowbreak 1+q+2q^{2}+3q^{3}+5q^{4}+7q^{5}+ \\
&&11q^{6}+15q^{7}+22q^{8}+\mathcal{O}\left( q^{9}\right). \end{eqnarray*}

When $g=0,$ our main theorem proves the formula predicted by Yau and Zaslow
\cite{Y-Z} for primitive classes.  For $g\geq 0$, G\"ottsche has recently
conjectured  a very general set of
formulas for the 
number of curves on algebraic surfaces \cite{Gott-conj} and Theorem \ref{thm:
main thm} proves his conjecture for primitive classes in  $K3$ surfaces.
Yau and Zaslow give a beautiful, though
indirect, argument that would be a
complete proof for the $g=0$ case if one could control the
complexity of singularities that can occur for curves in a complete linear
system. Beauville \cite{Beau}, Chen \cite{Chen}, and
Fantechi-G\"ottsche-van~Straten \cite{F-G-vS} have partial results
along this line.

We shall use a completely different argument by studying Gro\-mov-Wit\-ten
invariants for the twistor family of symplectic structures on a $K3$ surface.
We learned the twistor family approach from Li and Liu \cite{Li-Liu} who
studied the Seiberg-Witten theory for families and obtained interesting
results.

In the case of a hyperk\"ahler $K3$, the twistor family is the unit sphere
in the space of self-dual harmonic 2-forms. The idea of extending the
moduli space of pseudo-holomorphic curves by including the family of
non-degenerate, norm 1, self-dual, harmonic 2-forms goes back to Donaldson
\cite{Do}. He pointed out that in order to have the theory of
pseudo-holomorphic curves on a 4-manifold more closely mimic the theory of
divisors on a projective surface, one should include this family.

%In \cite{Kont}, Kontsevich also suggests using this family and
%attributes the idea to D. Morrison.

One key point in the proof of our main theorem is the use of the large
diffeomorphism group of a $K3$ surface to 
move $C$ to a particular class $S+\left( g+n\right) F$ on an elliptic $K3$
surface with section $S$ and fiber $F$ which has $24$ nodal fibers. Inside
the linear system $\left| S+\left( g+n\right) F\right| $, we can completely
understand the moduli space of stable maps and really
{\em count} the invariants $N_{g}(n)$, reducing the calculation to the
computation of ``local contributions'' by multiple covers.

The contribution of multiple covers  of the
smooth fibers is responsible for the $dG_{2}/dq$ term and the contribution
from multiple covers of the nodal fibers is related to the partition
function $p\left( d\right)$. The computation for the multiple covers of
nodal fibers requires an obstruction bundle computation. This is done by
splitting the moduli space into components and then identifying each
component with a moduli-obstruction problem arising from the Gromov-Witten
invariants of a certain blow-up of $\P ^{2}$. This ``matching'' technique
allows us to use known properties of the Gromov-Witten invariants of $\P
^{2}$ blown-up, specifically their invariance under Cremona
transformations, to show that the contribution of each component is always
0 or 
1. The partition function then arises combinatorially in a somewhat unusual
way (see lemma \ref{lem: no. of 1-admissable seqs is p(a)}). These
computations occupy section \ref{sec: analysis and local computations}.

These invariants $N_{g}(n)$ are all enumerative (Theorem \ref{thm: the
invariants are enumerative}). 

We can also apply our method to other elliptic surfaces. We blow up
$\mathbf{P}  ^{2}$ at nine distinct points and call the resulting algebraic
surface $Y$.  We consider the Gro\-mov-Wit\-ten invariant $N_{g}^{Y}(C)$ which
counts the number of curves of geometric genus $g$ passing
through $g$ generic points in a fixed class $C$. We show that any
class $C\in H_{2}(Y)$ 
whose genus $g$ invariants require exactly $g$ point constraints is related
to a 
class  of the form $C_{n}=\left( 
g+n\right) \left[ 3h-\sum_{i=1}^{9}e_{i}\right] +e_{9}$ by a Cremona
transform. Here $h$ is the
pullback to $Y$ of the hyperplane class in $\mathbf{P}^{2}$ and
$e_{1},...,e_{9}$ are the exceptional curves in $Y$. We obtain the
following:          

\begin{theorem}\label{thm: E1 case}
Let $Y$ be the rational elliptic surface; for fixed $g$ let $C_{n}$ be the
class 
$S+(g+n)F$ where $F$ is the fiber and $S$ is any section.
Let $N^{Y}_{g}(C_{n})$ denote the number of genus $g$ curves in the class
of $C_{n}$ passing through $g$ generic points. Then
\begin{eqnarray*}
\sum_{n=0}^{\infty }N_{g}^{Y}( C_{n}) q^{n}&=&\left(\sum_{k=1}^{\infty
}k(\sum_{d|k}d)q^{k-1} \right)^{g}\prod_{m=1}^{\infty
}\left( 1-q^{m}\right) ^{-12}\\
&=&\left(\frac{d}{dq}G_{2}(q) \right)^{g}\left(\frac{q}{\Delta (q)}
\right)^{1/2 }.
\end{eqnarray*}
\end{theorem}

Note that the Euler characteristic of $Y$ is 12 while the Euler
characteristic of $K3$ is 24. For us the relevant manifestation of this
fact is that elliptically fibered $K3$ surfaces have (generally) 24 nodal
fibers while rational elliptic surfaces have (generally) 12 nodal fibers. 

Since $N_{g}^{Y}( C_{n})$ is an ordinary Gro\-mov-Wit\-ten
invariant (without family), it is an invariant for the deformation class of
the symplectic structure on $Y$. In particular, $N_{g}^{Y}( C_{n})$
is independent of the locations of those blow up points in
$\mathbf{P}^{2}$ and it is left invariant by Cremona transforms. 
For the genus zero case, the invariants were obtained by G\"{o}ttsche and
Pandharipande \cite{Go-Pa} where they computed the
quantum cohomology for $\mathbf{P}^{2}$ blown up at arbitrary number of points
using the associativity law. Their numbers are in terms of two complicated
recursive formulas and it is not obvious that the numbers that correspond
to $N_{0}^{Y}( C_{n})$ can be put together to form modular forms, but Theorem
\ref{thm: E1 case} can be verified term by term for $g=0$ using their
recursion relations. 

The foundations on which our calculations rest have been developed by Li  
and Tian (\cite{Li-Tian3}\cite{Li-Tian2}\cite{Li-Tian}). They construct the
virtual fundamental cycle of the  moduli 
space of stable maps both symplectically and algebraically and they show that 
the two contructions coincide in the projective case. Ionel and Parker have a 
different approach to computing $N_0^Y(C_n)$ that does not rely on
\cite{Li-Tian3}.

Although our methods are completely different from those of Yau and Zaslow,
for the sake of completeness we ouline their beautiful argument 
for counting rational curves with $n$ nodes. Choose a smooth curve $C$ in
the $K3$ surface 
$  X$ with $C\cdot C=2n-2.$ By adjunction formula, the genus of $C$ equals
$n.$ One can show that $C$ moves in a complete linear system of dimension
$n$ using the Riemann-Roch theorem and a vanishing theorem. That is $\left|
C\right| \cong \mathbf{P}^{n}.$         

Imposing a node will put one constraint on the linear system $\left|
C\right| $. Therefore, by imposing $n$ nodes, one expects to see a finite
number of rational curves with $n$ nodes. Define $N_{0}(n)$ to be this
number. Now look at the compactified universal Jacobian $\pi
:\mathcal{\bar{J  }}\rightarrow \left| C\right| $ for this linear system
(c.f. Bershadsky, Sadov, and Vafa \cite{vafa}).
It is a smooth hyperk\"{a}hler manifold of dimension $2n.$          

If one assumes that each member in the linear system $\left| C\right| $ has
at most nodal singularities, then one can argue that for any $C^{\prime
}\in \left| C\right| $ the Euler characteristic of $\pi ^{-1}\left(
C^{\prime }\right) $ is always zero unless $C^{\prime }$ is a rational
curve with $n$ nodes. In the latter case, the Euler characteristic of $\pi
^{-1}\left( C^{\prime }\right) $ equals one. One concludes that $\chi
\left( \mathcal{  \bar{J}}\right) =N_{0}(n).$

One the other hand, $\mathcal{\bar{J}}$ is birational to the Hilbert scheme
$  \mathcal{H}_{n}$ of $n$ points in $X,$ which is again another smooth
hyperk\"{a}hler manifold. Using a result of  Batyrev
\cite{Bat} which
states that compact, birationally equivalent, projective, Calabi-Yau
manifolds have the same Betti numbers, one can conclude that $N_{0}(n)=\chi
\left(\mathcal{H}_{n}\right)$.  

Then one uses the result of G\"{o}ttsche \cite{Gott}, who
used Deligne's answer to the Weil conjecture to compute (among other
things) the Euler characteristic of $\mathcal{H}_{n}:$          
\[
\sum_{n=0}^{\infty }\chi \left( \mathcal{H}_{n}\right)
q^{n}=\prod_{m=1}^{\infty }\left( 1-q^{m}\right) ^{-24}. \] 
Finally, combining these results, one obtains
$$F_{0}=\sum_{n=0}^{\infty}N_{0}(n)q^{n}=\prod_{m=1}^{\infty }\left(
1-q^{m}\right) ^{-24}.$$

To turn this argument into a proof one must address the problems of 
what kind of singularities can occur in $\left| C\right| $ 
and how  the compactified Jacobian contributes to the Euler
characteristic of $\mathcal{\bar{J}}$.
Work on these issues has been started by
Beauville \cite{Beau}, Chen \cite{Chen}, and
Fantechi-G\"ottsche-van~Straten \cite{F-G-vS}.  

Our proof circumvents these problems by using symplectic geometry and
Gro\-mov-Wit\-ten invariants for families of symplectic structures. Our method
is also more 
direct than the Yau-Zaslow argument, avoiding the characteristic $p$
methods employed by
G\"ottsche and Batyrev.

We end this introduction with some speculation into the meaning of our
results. 
Ordinary Gro\-mov-Wit\-ten invariants give rise to the quantum cohomology ring.
There may be a corresponding structure in the
context of Gro\-mov-Wit\-ten invariants for families. The computation of
Theorem \ref{thm: main thm} shows that there is structure amongst these
invariants and suggests that there should be an interesting theory of
quantum cohomology that encodes it.

The ordinary quantum cohomology ring of $X$ gives a Frobenius structure on
$H^{*}(X;\mathbf{C}) $ and the (generalized) mirror conjecture states that
this Frobenius structure is equivalent to a Frobenius structure arising
from some sort of ``mirror object''  (see Givental \cite{Giv}). In the case
of a Calabi-Yau 3-fold, 
the mirror object is another Calabi-Yau 3-fold and the Frobenius structure
arises from its variation of Hodge structure. Theorem \ref{thm: main thm}
shows that the Gro\-mov-Wit\-ten invariants for $K3$ with its 
twistor family can be expressed in terms of quasi-modular forms. If there
is a quantum cohomology theory associated to the 
Gro\-mov-Wit\-ten invariants for families and a corresponding mirror
conjecture, then our theorem should provide clues as to what the ``mirror
object'' of $K3$ with its twistor family should 
be.

This paper is organized as follows. In section \ref{sec: inv of families}
we define invariants for families of symplectic structures; in section
\ref{sec: twistor families} we discuss twistor families and define
$N_{g}(n)$; in section \ref{sec: computation of Ngn} we compute $N_{g}(n)$
and prove our main theorem; in section
\ref{sec: analysis and local computations} we analyze the moduli spaces and
compute local contributions; in section \ref{sec: counting on E1} we apply
similar techniques to $\mathbf{P}^{2}$ blown up at nine points.

The authors are pleased to acknowledge helpful conversations with A.
Bertram, A. Givental, L. G\"ottsche, E. Ionel, A. Liu, P. Lu, D. Maclagan,
D. McKinnon, T. Parker, S. Schleimer, C. Taubes, A. Todorov, 
and S.-T. Yau.  The authors especially thank L. G\"ottsche for
sharing early versions of  his conjecture with us and for providing many
other valuable 
communications. We would like to thank R. Pandharipande and L. G\"ottsche
for sending us their Maple program to verify our results. Additionally 
we  thank the Park City Mathematics Institute for support and providing a
stimulating environment where part of this work was carried out.

\section{Invariants of families of symplectic structures}\label{sec: inv of
families} 

In this section, we introduce an invariant for a family of symplectic
structures $\omega _{B}:B\rightarrow \Omega _{sympl}^{2}\left( X\right) $
on a compact manifold $X.$ Here $B$ is an oriented, compact manifold and
$\omega _{B}$      
is a smooth map into the space of symplectic forms $\Omega
_{sympl}^{2}\left( X\right) .$ This invariant is a direct generalization of
the Gro\-mov-Wit\-ten invariants. Roughly speaking, it counts the number of
maps $u:\Sigma \to X$ which are holomorphic with respect to some almost
complex structure in a generic family compatible with $\omega _{B}.$
Kronheimer \cite{Kr} and Li and Liu \cite{Li-Liu} have also
studied invariants for families of symplectic structures and obtained
interesting results. 

In their paper on Gro\-mov-Wit\-ten invariants for general symplectic
manifolds, Li and Tian \cite{Li-Tian} setup a general framework for
constructing invariants. Their results are easy to adapt to our setting but
we remark that in the case of interest, the full Li-Tian machinery is not
needed and the techniques of Ruan-Tian 
\cite{Ruan-Tian} would suffice. This is because 2-dimension families of
symplectic structures on a 4-manifold behave like the semi-positive case
for the ordinary invariants, \ie the (perturbed) moduli spaces are
compactified by strata of codimension at least 2. The definition of
Gro\-mov-Wit\-ten invariants for families is also contained as part of the very
general approach of Ruan \cite{Ruan}. We employ the Li-Tian machinery
because they are also able to relate their symplectic constructions to
their purely algebraic ones (\cite{Li-Tian2} and \cite{Li-Tian3}). 

Let $X$ be a compact smooth manifold. Suppose that $\omega _{B}$ is a
smooth family of symplectic structures on $X$ parameterized by an oriented,
compact manifold $  B. $ Let $J_{B}:B\rightarrow \mathcal{J}\left( X\right)
$ be a smooth family of almost complex structures on $X$ such that
$J_{t}=J_{B}\left( t\right) $ is compatible with $\omega _{t}=\omega
_{B}\left( t\right) $ for any $t\in B.$ In particular, $g_{t}=\omega
_{t}\left( \cdot ,J_{t}\cdot \right) $ is a family of Riemannian metrics on
$X.$ It is not difficult to see that $J_{B}$ always exists and is unique up to
homotopy. This follows from the fact that the space of all almost complex
structures compatible with a fixed symplectic form is contractible.         

Given $X$ and $\omega _{B}$ as above, we shall define the GW-invariant for
family as a homomorphism: 
\[
\Psi _{\left( A,g,k\right) }^{(X,\omega
_{B})}:\bigoplus_{i=1}^{k}H^{a_{i}}\left( X,\mathbf{Q}\right) \bigoplus
H^{b}\left( \overline{\mathcal{M}}_{g,k},\mathbf{Q}\right) \rightarrow \mathbf{Q},
\]         
with 
\begin{equation}\label{eqn: dimension formula}
\sum_{i=1}^{k}a_{i}+b=2c_{1}\left( X\right) \left( A\right) +2k+\dim
B+\left( \dim X -6\right) \left( 1-g\right) .
\end{equation}

 Here $A\in H_{2}\left( X,\mathbf{Z}
\right) $ and $\overline{\mathcal{M}}_{g,k}$ is the Deligne-Mumford
compactification of the moduli space of Riemann surfaces of genus $g$ with
$  k $ distinct marked points (define $\overline{\mathcal{M}}_{g,k}$ to be
a point if $2g+k<3$).         

For any particular symplectic structure $\omega _{t}$ and the corresponding
almost complex structure $J_{t},$ Li and Tian define a section $\Phi _{t}$
of $E\rightarrow \overline{\mathcal{F}_{A}}$ $\left( X,g,k\right) $ (we
recall the definition below) which
is equivalent to the Cauchy-Riemann operator. These sections depend on
$t\in B$ smoothly so that we have a section $\Phi $ of $E\rightarrow
\overline{  \mathcal{F}}_{A}$ $\left( X,g,k\right) $ $\times B.$          

Let us first recall their notations:

A stable map with $k$ marked points is a tuple $\left( f,\Sigma
;x_{1},...x_{k}\right) $ satisfying:
\begin{enumerate}
\item [{\bf (i)}] $\Sigma
=\bigcup\limits_{i=1}^{m}\Sigma _{i}$ is a connected normal crossing
projective curve and $x_{i}$'s are distinct smooth points on $\Sigma$,
\item [{\bf (ii)}] $f$ is continuous and $f|_{\Sigma _{i}}$ can be lifted to a
smooth map on the normalization of $\Sigma _{i}$, and
\item [{\bf (iii)}] if $\Sigma _{i}$ is a
smooth rational curve such that $f\left( \Sigma _{i}\right) $ represents a
trivial homology class in $H_{2}\left( X,\mathbf{Q}\right) ,$ then the
cardinality of $  \Sigma _{i}\bigcap \left( \left\{
x_{1},...,x_{k}\right\} \bigcup S\left( \Sigma \right) \right) $ is at
least three where $S\left( \Sigma \right) $ is the singular set of
$\Sigma.$    
\end{enumerate}

Two stable maps $\left( f,\Sigma ;x_{1},...x_{k}\right) $ and $\left(
f^{\prime },\Sigma ;x_{1}^{\prime },...x_{k}^{\prime }\right) $ are
equivalent if there is a biholomorphism $\sigma :\Sigma \rightarrow \Sigma
^{\prime }$ such that $\sigma \left( x_{i}\right) =x_{i}^{\prime }$ for $
1\leq i\leq k$ and $f^{\prime }=f\circ \sigma .$ We denote the space of
equivalent classes of stable maps of genus $g$ with $k$ marked points and
with total homology class $A$ by $\overline{\mathcal{F}}_{A}$ $\left(
X,g,k\right) $ and the subspace consisting of equivalent classes of stable
maps with smooth domain by $\mathcal{F}_{A}\left( X,g,k\right) .$ The
topology of $\overline{\mathcal{F}}_{A}$ $\left( X,g,k\right) $ can be
defined by sequential convergence. Next they introduced a generalized
bundle $E$ over $\overline{\mathcal{F}}_{A}$ $\left( X,g,k\right) $ as
follows: For any $\left[ \left( f,\Sigma ;x_{1},...x_{k}\right) \right] \in
\overline{  \mathcal{F}}_{A}$ $\left( X,g,k\right) ,$ the fiber of $E$
consists of all $  f^{*}TX$-valued $\left( 0,1\right) $-forms over the
normalization of $\Sigma . $ Equiped with the continuous topology, $E$ is
a generalized bundle over $  \overline{\mathcal{F}}_{A}$ $\left(
X,g,k\right) $ in the sense of Li and Tian.                  

For each $t\in B,$ there is a section of $E$ given by the Cauchy-Riemann
operator defined by $J_{t}:$ Namely, for any $\left[ \left( f,\Sigma
;x_{1},...,x_{k}\right) \right] \in \overline{\mathcal{F}}_{A}$ $\left(
X,g,k\right) ,$ we have $\Phi _{t}\left( f,\Sigma ;x_{1},...x_{k}\right)
=df+J_{t}\circ df\circ j_{\Sigma }$ where $j_{\Sigma }$ is the complex
structure of $\Sigma .$ Putting different $t\in B$ together, we have a
section $\Phi $ of $E$ over $\overline{\mathcal{F}}_{A}$
$\left(X,g,k\right) \times B$ given by         
\[
\Phi \left( \left[ \left( f,\Sigma ;x_{1},...x_{k}\right) \right] ,t\right)
=df+J_{t}\circ df\circ j_{\Sigma }. 
\]

The following theorems are easy adaptations of those in of Li and Tian
found in 
\cite{Li-Tian}:

\begin{theorem}
The section $\Phi $ gives rise to a generalized Fredholm orbifold bundle
with the natural orientation and of index 
$$2c_{1}( X) \left[
A\right] +2k+\dim B+\left( \dim X -6\right) \left( 1-g\right). $$
\end{theorem}

\begin{theorem}
Let $\omega _{B}$ and $\omega _{B}^{\prime }$ are two families of
symplectic structures on $X$ parameterized by $B$. Suppose that they are
equivalent to each other under deformations for families. Let $J_{B}$ and
$J_{B}^{\prime }$ be two families of almost complex structures on $X$
compatible with corresponding symplectic structures.         

Suppose that $\Phi $ and $\Phi ^{\prime }$ are the corresponding section of
$  E$ over $\overline{\mathcal{F}}_{A}$ $\left( X,g,k\right) \times B.$
Then $  \Phi $ and $\Phi ^{\prime }$ are homotopic to each other as
generalized Fredholm orbifold bundles.         
\end{theorem}

Using the main theorem of Li and Tian in their paper, there is an Euler
class $e\left( \left[ \Phi :\overline{\mathcal{F}}_{A}\left( X,g,k\right)
\times B\rightarrow E\right] \right) $ in $H_{r}\left(
\overline{\mathcal{F}  _{A}}\left( X,g,k\right) \times B,\mathbf{Q}\right) $
with $r=2c_{1}\left( X\right) \left[ A\right] +2k+\dim B+\left( 
\dim X -6\right)
\left( 1-g\right)$.  This class is called the {\em virtual fundamental cycle} 
of the moduli space of holomorphic stable maps $\mathcal{M}=\M (X,\omega
_{B},g,k,C)$.  We denote
it by $[\mathcal{M}]^{vir}$.

To define the invariant for $\omega _{B},$ we consider the following two
maps. First we have the evaluation map $ev:\overline{\mathcal{F}}_{A}\left(
X,g,k\right) \times B\rightarrow X^{k}$ : \[         
ev\left( \left( f,\Sigma ;x_{1},...,x_{k}\right) ,t\right) =\left( f\left(
x_{1}\right) ,...,f\left( x_{k}\right) \right) . \] 
and second we have the forgetful map $\pi _{g,k}:\overline{\mathcal{F}}_{A}
\left( X,g,k\right) \times B\rightarrow \overline{\mathcal{M}}_{g,k}:$  \[ 
\pi _{g,k}\left( \left( f,\Sigma ;x_{1},...,x_{k}\right) ,t\right)
=red\left( \Sigma ;x_{1},...,x_{k}\right) . \] 
Here $red\left( \Sigma ;x_{1},...,x_{k}\right) $ is the stable reduction of
$  \left( \Sigma ;x_{1},...,x_{k}\right) $ that is obtained by contracting
all of its non-stable irreducible components.  

Now we can define the invariants
\[
\Psi _{\left( A,g,k\right) }^{(X,\omega _{B})}:H^{*}\left( X,\mathbf{Q}\right)
^{k}\times H^{*}\left( \overline{\mathcal{M}}_{g,k},\mathbf{Q}\right)
\rightarrow 
\mathbf{Q}.  
\]
by
\[
\Psi _{\left( A,g,k\right) }^{(X,\omega _{B})}\left(\alpha
_{1},...,\alpha _{k};\beta \right) =\left( ev^{*}\left( \pi _{1}^{*}\alpha
_{1}\wedge ...\wedge \pi _{k}^{*}\alpha _{k}\right) \cup \pi
_{g,k}^{*}\left( \beta \right) \right) [\mathcal{M}]^{vir} \]
for any $\alpha _{1},...,\alpha _{k}\in H^{*}\left( X,\mathbf{Q}\right) $ and
$  \beta \in H^{*}\left( \overline{\mathcal{M}}_{g,k},\mathbf{Q}\right) .$ 

\begin{theorem}
$\Psi _{\left( A,g,k\right) }^{(X,\omega _{B})}$ is an invariant of the
deformation class of the family of symplectic structures $\omega _{B}.$
\end{theorem}   

If $(\hat{\alpha }_{1},\ldots,\hat{\alpha }_{k})$ are geometric cycles in
$X$ that are Poincar\'e dual to $(\alpha _{1},\ldots,\alpha _{k})$ and
$\hat{\beta }$ is a cycle in $\overline{\mathcal{M}}_{g,k}$ dual to $\beta
$, then $\Psi ^{(X,\omega _{B})}_{(A,g,k)}(\alpha _{1},\ldots,\alpha
_{k};\beta )$ counts the number of stable maps $f:\Sigma _{g}\to X$ so that \begin{enumerate}
\item $f(\Sigma _{g})$ represents the class $A$,
\item $f$ is $J_{t}$-holomorphic for some $t\in B$, 
\item $f(x_{i})$ lies on $\hat{\alpha }_{i}$, and
\item the stable reduction of $\Sigma _{g}$ lies in $\hat{\beta }\subset
\overline{\mathcal{M}}_{g,k}$.
\end{enumerate}

One is usually interested in $\hat{\beta }=\overline{\mathcal{M}}_{g,k}$,
or sometimes $\hat{\beta }=\point \in \overline{\mathcal{M}}_{g,k}$.

\section{Twistor families of $K3$ surfaces and the definition of
$N_{g}(n)$}\label{sec: twistor families}

In this section we collect some general facts about $K3$ surfaces and their
twistor families. We show that every twistor family is deformation
equivalent and we define $N_{g}(n)$ in terms of the Gro\-mov-Wit\-ten invariant
for this family. We show that when $X$ is projective and $|C|$ has only
reduced and irreducible curves, $N_{g}(n)$ coincides with
the enumerative count defined by algebraic geometers (see
\cite{F-G-vS}). The results of the 
section are summarized in Definition \ref{def: definition of Ng,n}.

A $K3$ surface is a simply-connected, compact, complex surface $X$ with
$c_{1}(X)=0$. For a general reference on $K3$ surfaces and twistor families
we refer the reader to \cite{BPV} or \cite{Bess}. Any pair of $K3$ surfaces
are deformation equivalent and 
hence diffeomorphic. A {\em marking} of a $K3$ surface $X$ is an
identification of $(H^{2}(X;\znums),Q_{X})$ with the fixed unimodular form
$Q=-2E_{8}\oplus 3\left(\begin{smallmatrix}
0&1\\
1&0
\end{smallmatrix} \right)$.  The space of marked $K3$ surfaces forms a
connected, 20 
complex dimensional moduli space.

The complex structure on a marked $K3$ surface $X$ is
determined by how the line $H^{0,2}(X)$ sits in $Q\otimes \cnums $. To make
this precise, define the {\em period} $\Omega _{X}$ of $X$ to be the
element of the {\em period domain} 
$$\mathcal{D} =\{\Omega \in \P (Q\otimes
\cnums ): \left\langle \Omega ,\overline{\Omega }\right\rangle
=0,\left\langle \Omega ,\Omega \right\rangle>0 \}$$
given by the image of $H^{0,2}(X)$ under the marking.

The Torelli theorem states that $\mathcal{D}$ is the moduli space of marked
$K3$ surfaces, \ie every marked $K3$ surface corresponds uniquely to its
period point in $\mathcal{D}$ and every $\Omega \in \mathcal{D}$ is the
period point of some $K3$ surface.

For a fixed $K3$ surface $X$ with period  $\Omega _{X} $, a class $\omega
_{X}\in Q\otimes \cnums $ is a K\"ahler 
class for $(X,\Omega _{X})$ if and only if $\left\langle \omega _{X},\Omega
_{X}\right\rangle=0$, $\left\langle \omega
_{X},\overline{\Omega}_{X}\right\rangle=0$, $\omega
_{X}=\overline{\omega}_{X}$, and   
$\left\langle \omega _{X},\omega _{X}\right\rangle>0$. For any K\"ahler
$K3$ surface $(X,\Omega _{X},\omega _{X})$ there is a unique hyperk\"ahler
metric by Yau's proof of the Calabi conjecture \cite{Yau}. A hyperk\"ahler
metric $g$ 
determines a 2-sphere worth of K\"ahler structures, namely the unit sphere
in the space $\mathcal{H}^{2}_{+,g}$ of self-dual harmonic forms. We can
describe the corresponding 2-sphere of period points as follows. Consider
the projective plane spanned by
$\left\langle \Omega _{X},\overline{\Omega }_{X}, \omega _{X}\right\rangle$
in $\P (Q\otimes \cnums) $. Since $\left\langle \Omega
_{X},\overline{\Omega }_{X}, \omega _{X}\right\rangle$ spans
$\mathcal{H}^{2}_{+,g}\otimes \cnums$, the intersection of this projective
plane and the period domain is the quadric determined by $\left\langle
\Omega ,\overline{\Omega }\right\rangle=0$. This is a smooth plane quadric
and hence a 2-sphere. This 2-sphere of complex structures together with the
corresponding 2-sphere of K\"ahler structures we call a {\em twistor
family}. We will use the notations $(J_{T},\omega_{T})$ to refer to a twistor
family and  $J_{t}$, $\omega_{t}$ for $t\in T$ to refer to individual
members. 

The following proposition was explained to us by Andrei Todorov:

\begin{prop}\label{prop: one curve in class of C in twistor family}
Let $(X,\Omega _{X},\omega _{X})$ be a marked, K\"ahler $K3$ surface and
$(J_{T},\omega_{T})$ the corresponding twistor family. Let $C\in
H^{2}(X;\znums  )$ be a class of square $C^{2}\geq -2$. Then there is
exactly one member $t\in T$ for which there is a $J_{t}$-holomorphic curve
in the class of $C$.
\end{prop}
{\em Proof:} A class $C$  of square $-2$ or larger admits a holomorphic
curve if and only if $C\in H^{1,1}(X;\znums )$ and $C$ pairs positively
with the 
K\"ahler class. Since $C$ is a real class, $C\in H^{1,1}$ if and
only if $\left\langle C,\Omega_{X}\right\rangle=0$. This equation determines a
hyperplane in $\P (H^{2}_{+})$ and so meets the twistor space in 2 points
$\pm \Omega_{0}$ (since the twistor space is a quadric). Then exactly one of
$\pm\Omega_{0}$ will have its corresponding K\"ahler class pair positively
with $C$.\qed 

We next show that every twistor family is the same up to deformation.

\begin{proposition}\label{prop: every twistor family is homotopic}
Let $X_{1}$ and $X_{2}$ be two K\"ahler $K3$ surfaces, then the
corresponding twistor families $T_{0}$ and $T_{1}$ are deformation
equivalent. 
\end{proposition}
{\em Proof:}  
The moduli space of $K3$ surfaces is connected and the space of
hyperk\"ahler structures for a fixed $K3$ surface is contractible (it is
the K\"ahler cone). Therefore, the space parameterizing hyperk\"ahler $K3$
surfaces $(X,\Omega_{X},\omega _{X})$ is also connected. We can 
thus find a path $(X_{s},\Omega _{s},\omega _{s})$, $s\in [0,1]$ connecting
$X_{0}$ to $X_{1}$ where the twistor family of $\omega _{i}$ is $T_{i}$ for
$i=0,1$. By then associating to each hyperk\"ahler structure $\omega _{s}$
its twistor family $T_{s}$, we obtain a continuous deformation of $T_{0}$ to
$T_{1}$. \qed 

From this proposition we see that the Gro\-mov-Wit\-ten invariants for a
twistor family are independent of the choice of a twistor family. We can
thus write unambiguously
$$\Psi  ^{(K3,\omega _{T})}_{(C,g,k)}:H^{*}(K3;\znums )^{k}\otimes
H^{*}(\overline{\mathcal{M}}_{g,k})\to \qnums  .$$

We are primarily interested in the invariants that count stable maps
without fixing the complex structure on the domain. That is, the invariants
obtained using the Poincar\'e dual of the fundamental class of
$\overline{\mathcal{M}}_{g,k}$. It is enough to consider those constraints
that come from the generator of $H^{4}(K3,\znums )$; these count curves
passing 
through fixed generic points. The invariants with the constraint that
the $k$th point lies on a fixed generic cycle dual to an element $\beta \in
H^{2}(K3)$ can be computed in terms of the invariants for $k-1$ constraints
and the pairing $\beta \cdot C$. For this reason, constraining the
invariants by elements of $H^{2}(K3)$ is uninteresting, and of course
elements of $H^{0}(K3)$ provide no constraints at all. 

We can thus simplify notation by defining
$$\Psi(C,g,k)\equiv \Psi ^{(K3,\omega
_{T})}_{(C,g,k)}(PD(x_{1}),\ldots,PD(x_{k});PD([\overline{\mathcal{M}}_{g,k}
])). $$ 

An important observation about the twistor family is the following.
\begin{prop}\label{prop: diffeos preserve T}
If $f:K3\to K3$ is an orientation preserving diffeomorphism, then
the pullback family $f^{*}(\omega _{T})$ is deformation equivalent to
$\omega _{T}$; thus 
$$
\Psi (C,g,k) =\Psi (f_{*}(C),g,k).  
$$
\end{prop}
{\em Proof:} Let $\omega _{T}$ be the twistor family associated to a
hyperk\"ahler metric $g$. Then $f^{*}(\omega _{T})$ is the twistor family
associated to the hyperk\"ahler metric $f^{*}(g)$ and so Proposition
\ref{prop: every twistor family is homotopic} they are deformation
equivalent. \qed 

The $K3$ surface has a big diffeomorphism group in the sense of Friedman and
Morgan \cite{Fr-Mo}, which means that every automorphism of the lattice $Q_X$
which preserves spinor norm can be realized by an orientation preserving
diffeomorphism. In particular, one can take any primitive class $C\in
H_2(K3;\znums)$ to any other primitive class with the same square via an
orientation preserving
diffeomorphism. 

We are now in a position to define $N_{g}(n)$. By the adjunction formula, a
holomorphic curve of genus $g$ with $n$ nodes will be in a class $C$ with
square $C^{2}=2(g+n)-2$.

\begin{definition}\label{def: definition of Ng,n}
 Let $C$ be any primitive class with
$C^{2}=2(g+n)-2$. We define the number $N_{g}(n)$ by
\begin{eqnarray*}
N_{g}(n)&=&\Psi (C,g,g)\\
&=&\Psi ^{(K3,\omega_T)}_{(C,g,g)}
(PD(x_{1}),\ldots,PD(x_{g});PD([\overline{\mathcal{M}}_{g,k}])).    
\end{eqnarray*}
By Proposition \ref{prop: diffeos preserve T} , $N_{g}(n)$ is independent of
the choice of the primitive class $C$. By Proposition \ref{prop: every
twistor family is homotopic}, $N_{g}(n)$ is independent of the choice of
twistor family. Finally, in the case of a projective $K3$ surface with an
effective divisor in the class of $C$, Proposition \ref{prop: one curve in
 class of C in 
twistor family} shows that $N_{g}(n)$ counts holomorphic maps $f:D\to X$ of
genus $g$ curves to $X$ with image in $|C|$ and passing through $g$
 generic points.  
\end{definition}

\begin{thm}\label{thm: the invariants are enumerative}
If $X$ is generic among those $K3$ surfaces admitting a curve in the class
$[C]$, then invariant $N_{g}(n)$ is enumerative.
\end{thm}
{\sc Proof:}\footnote{This argument is due to
Lothar G\"ottsche. We are grateful to him for showing it to us.}
 The assumption that the $K3$ surface
$X$ is generic among those admitting a curve in the class of $C$ guarantees
that
the primitive class $C$ generates the Picard group. 
Suppose that the invariant $\Psi
(C,g,g)$ differs from the actual count of curves $\Sigma \in |C|$ of genus
$g$ 
passing through $g$ general points (curves counted with multiplicities if
the they are not all nodal as in \cite{F-G-vS}).  Then there 
is an extra map $D\to X$ of a curve of arithmetic genus $g$ that has some
contracted components and the rest of the
map is generically injective with irreducible image. Let $C_{1}$ be a
contracted component. Since the $g$ marked points have to go to $g$
distinct points on $X$, $C_{1}$ can have at most 1 marked point. By
stability then, either the geometric genus $g(C_{1})$ is larger than 0 or
$C_{1}$ intersects the rest of $D$ in at least 2 points. Since the image of
$D$ is irreducible, the contracted components cannot all be genus 0
unless the dual graph of $D$ is not a tree. Thus either $D$ has a
contracted component of genus greater than 0 or the dual graph of $D$ is
not a tree. In either case, the geometric genus of the image is smaller
than the arithmetic genus of $D$ and thus the image is a curve of genus
less than $g$ passing through $g$ points. This does not occur for $g$
generic points by a dimension count.

\begin{rem}\label{rem: our method only is for primitive C}
The conjectured formula of Yau and Zaslow applies to non-primitive classes
as well. The above definition could be made for arbitrary classes $C$, but
{\it a priori} $N_{g}(n)$ would also depend on the divisibility of $C$. Our
method of computing $N_{g}(n)$ only applies to primitive classes, so the
Yau-Zaslow conjecture remains open for the non-primitive classes. 
\end{rem}

\section{Computation of $N_{g}(n)$}\label{sec: computation of Ngn}

To compute $N_{g}(n)$ we are free to choose any family of symplectic
structures deformation equivalent to the twistor family and any
primitive class $C$ with $C^{2}=2(g+n)-2$. Let $X$ be
an elliptically fibered $K3$ surface with a section and 24 nodal singular
fibers $N_{1},\ldots,N_{24}$. Endow $X$ with a hyperk\"ahler metric and let
$(\omega _{T},J_{T})$ be the 
corresponding twistor family. Let $S$ denote the section and $F$ the class of
the fiber so that $F^{2}=0$, $F\cdot S=1$, and $S^{2} =-2$. Let $C$ be the
class $S+(n+g)F$ and fix $g$ generic points $x_{1},\ldots,x_{g}$ not on $S$
that lie
on $g$ distinct smooth fibers which we label $F_{1},\ldots, F_{g}$ (see the
following illustration).

\vskip 15pt
\setlength{\unitlength}{0.00083333in}
\begingroup\makeatletter\ifx\SetFigFont\undefined
% extract first six characters in \fmtname
\def\x#1#2#3#4#5#6#7\relax{\def\x{#1#2#3#4#5#6}}%
\expandafter\x\fmtname xxxxxx\relax \def\y{splain}%
\ifx\x\y   % LaTeX or SliTeX?
\gdef\SetFigFont#1#2#3{%
  \ifnum #1<17\tiny\else \ifnum #1<20\small\else
  \ifnum #1<24\normalsize\else \ifnum #1<29\large\else
  \ifnum #1<34\Large\else \ifnum #1<41\LARGE\else
     \huge\fi\fi\fi\fi\fi\fi
  \csname #3\endcsname}%
\else
\gdef\SetFigFont#1#2#3{\begingroup
  \count@#1\relax \ifnum 25<\count@\count@25\fi
  \def\x{\endgroup\@setsize\SetFigFont{#2pt}}%
  \expandafter\x
    \csname \romannumeral\the\count@ pt\expandafter\endcsname
    \csname @\romannumeral\the\count@ pt\endcsname
  \csname #3\endcsname}%
\fi
\fi\endgroup
{\renewcommand{\dashlinestretch}{30}
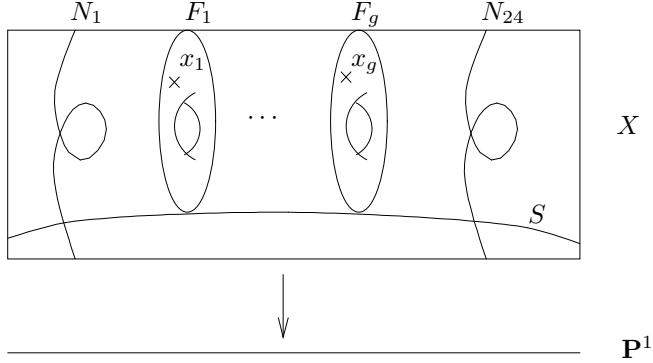
\begin{picture}(3928,2183)(0,-10)
\put(1286.942,1456.000){\arc{473.107}{2.0399}{4.2433}}
\put(1029.276,1439.500){\arc{367.450}{5.1978}{7.3686}}
\put(1131,1488){\ellipse{356}{1136}}
\put(2357.942,1456.000){\arc{473.107}{2.0399}{4.2433}}
\put(2100.275,1439.500){\arc{367.451}{5.1978}{7.3685}}
\put(2202,1488){\ellipse{356}{1136}}
\path(2093,1796)(2093,1796)(2154,1732)
\path(2093,1732)(2154,1796)
\path(1022,1764)(1022,1764)(1083,1699)
\path(1022,1699)(1083,1764)
\path(12,2056)(3582,2056)(3582,628)
	(12,628)(12,2056)
\path(12,43)(3582,43)
\path(1732,530)(1732,141)
\path(1702.000,261.000)(1732.000,141.000)(1762.000,261.000)
\path(12,758)	(83.423,779.758)
	(145.098,798.203)
	(198.071,813.615)
	(243.386,826.274)
	(315.231,844.458)
	(369.000,855.000)

\path(369,855)	(433.726,864.056)
	(507.440,872.363)
	(547.266,876.242)
	(588.859,879.939)
	(632.057,883.458)
	(676.700,886.801)
	(722.629,889.970)
	(769.682,892.968)
	(817.700,895.796)
	(866.522,898.457)
	(915.989,900.952)
	(965.938,903.286)
	(1016.211,905.458)
	(1066.648,907.472)
	(1117.086,909.331)
	(1167.368,911.035)
	(1217.331,912.588)
	(1266.817,913.992)
	(1315.664,915.248)
	(1363.712,916.360)
	(1410.801,917.329)
	(1456.771,918.157)
	(1501.461,918.848)
	(1544.712,919.402)
	(1586.362,919.823)
	(1626.252,920.113)
	(1700.109,920.306)
	(1765.000,920.000)

\path(1765,920)	(1834.318,919.048)
	(1872.665,918.313)
	(1913.243,917.405)
	(1955.879,916.325)
	(2000.401,915.073)
	(2046.638,913.649)
	(2094.417,912.054)
	(2143.566,910.287)
	(2193.915,908.349)
	(2245.290,906.241)
	(2297.521,903.963)
	(2350.434,901.514)
	(2403.859,898.895)
	(2457.623,896.107)
	(2511.555,893.149)
	(2565.482,890.022)
	(2619.234,886.726)
	(2672.637,883.262)
	(2725.520,879.629)
	(2777.712,875.828)
	(2829.039,871.859)
	(2879.331,867.723)
	(2928.416,863.420)
	(2976.121,858.949)
	(3022.275,854.312)
	(3066.706,849.508)
	(3109.243,844.538)
	(3149.712,839.402)
	(3187.942,834.100)
	(3257.000,823.000)

\path(3257,823)	(3306.288,811.936)
	(3371.785,793.474)
	(3413.002,780.738)
	(3461.140,765.275)
	(3517.154,746.793)
	(3582.000,725.000)

\path(434,2056)	(415.272,2011.008)
	(398.182,1969.006)
	(382.680,1929.827)
	(368.720,1893.305)
	(345.231,1827.571)
	(327.334,1770.477)
	(314.644,1720.700)
	(306.781,1676.912)
	(304.000,1602.000)

\path(304,1602)	(308.522,1557.794)
	(316.624,1505.002)
	(328.763,1447.856)
	(345.396,1390.590)
	(366.980,1337.438)
	(393.970,1292.633)
	(466.000,1245.000)

\path(466,1245)	(522.439,1258.463)
	(574.714,1304.815)
	(613.382,1369.760)
	(629.000,1439.000)

\path(629,1439)	(617.973,1496.038)
	(588.484,1549.574)
	(546.753,1588.573)
	(499.000,1602.000)

\path(499,1602)	(457.232,1589.703)
	(420.852,1561.611)
	(389.651,1521.575)
	(363.419,1473.445)
	(341.947,1421.072)
	(325.026,1368.307)
	(312.447,1318.999)
	(304.000,1277.000)

\path(304,1277)	(298.252,1226.573)
	(298.024,1171.174)
	(303.699,1108.903)
	(315.660,1037.861)
	(324.117,998.458)
	(334.288,956.149)
	(346.222,910.699)
	(359.966,861.868)
	(375.568,809.420)
	(393.076,753.118)
	(412.537,692.724)
	(434.000,628.000)

\path(2998,2056)	(2979.257,2011.008)
	(2962.152,1969.006)
	(2946.639,1929.827)
	(2932.669,1893.305)
	(2909.168,1827.571)
	(2891.266,1770.477)
	(2878.581,1720.700)
	(2870.730,1676.912)
	(2868.000,1602.000)

\path(2868,1602)	(2872.492,1557.794)
	(2880.573,1505.002)
	(2892.700,1447.856)
	(2909.329,1390.590)
	(2930.917,1337.438)
	(2957.920,1292.633)
	(3030.000,1245.000)

\path(3030,1245)	(3086.388,1258.463)
	(3138.646,1304.815)
	(3177.331,1369.760)
	(3193.000,1439.000)

\path(3193,1439)	(3181.922,1496.038)
	(3152.416,1549.574)
	(3110.702,1588.573)
	(3063.000,1602.000)

\path(3063,1602)	(3021.202,1589.703)
	(2984.800,1561.611)
	(2953.585,1521.575)
	(2927.347,1473.445)
	(2905.879,1421.072)
	(2888.971,1368.307)
	(2876.414,1318.999)
	(2868.000,1277.000)

\path(2868,1277)	(2862.222,1226.573)
	(2861.973,1171.174)
	(2867.636,1108.903)
	(2879.592,1037.861)
	(2888.050,998.458)
	(2898.225,956.149)
	(2910.164,910.699)
	(2923.915,861.868)
	(2939.527,809.420)
	(2957.046,753.118)
	(2976.521,692.724)
	(2998.000,628.000)

\put(1083,1829){\makebox(0,0)[lb]{\smash{$x_1$}}}
\put(2154,1829){\makebox(0,0)[lb]{\smash{$x_g$}}}
\put(401,2121){\makebox(0,0)[lb]{\smash{$N_1$}}}
\put(1115,2121){\makebox(0,0)[lb]{\smash{$F_1$}}}
\put(2154,2121){\makebox(0,0)[lb]{\smash{$F_g$}}}
\put(2965,2121){\makebox(0,0)[lb]{\smash{$N_{24}$}}}
\put(3809,1407){\makebox(0,0)[lb]{\smash{$X$}}}
\put(3842,11){\makebox(0,0)[lb]{\smash{$\mathbf{P}^1$}}}
\put(1505,1472){\makebox(0,0)[lb]{\smash{$\cdots$}}}
\put(3257,855){\makebox(0,0)[lb]{\smash{$S$}}}
\end{picture}
}
\vskip 15pt

Recall
from Definition \ref{def: definition of Ng,n}  that 
$$N_{g}(n)=\Psi (C,g,g). $$
This counts the number of stable maps of geometric genus $g$,
in the class 
of $C$, whose image contains $x_{1},\ldots,x_{g}$, and which are
holomorphic for some $J_{t}$, $t\in T$. The space of such maps forms a
compact moduli space of virtual dimension 0 which we denote by
$\mathcal{M}_{C,g}$. By Proposition \ref{prop: one curve in class of C in
twistor family}, there is a unique $t\in T$ so that there
are $J_{t}$-holomorphic curves in the class of $C$. This $J_{t}$ must be
the original 
elliptically fibered complex structure. Thus $\mathcal{M}_{C,g}$ consists
of stable holomorphic maps whose image are in the linear system
$|S+(n+g)F|$ and contain the points $x_{1},\ldots,x_{g}$. Because of the
elliptic fibration, the linear system $|C|$ is easy to analyze. 

The dimension of $|S+(n+g)F|$ is $n+g$ and consists solely of reducible
curves which are each a union of the section and $(n+g)$ (not necessarily
distinct) fibers. Since the image contains the points $x_{1},\ldots,x_{g}$,
it contains the corresponding smooth fibers $F_{1},\ldots,F_{g}$. The image
of a map in 
$\mathcal{M}_{C,g}$ must therefore be the union of the section $S$, the $g$
fibers $F_{1},\ldots,F_{g}$, and some number of  nodal fibers (possibly
counted with 
multiplicity). We summarize this discussion in the following

\begin{prop}
Let $\mathcal{M}_{C,g}$ be the moduli space of stable maps of genus $g$ in
the class of $C=S+(g+n)F$ passing through the points
$x_{1},\ldots,x_{g}$. Let $\pi :\mathcal{M}_{C,g}\to \mathbf{P}(H^{0}(X,C))
$ be the natural projection onto the linear system $|C|$. Then $\im (\pi )$
is a finite number of points labeled by the vectors
$\mathbf{a}=(a_{1},\ldots,a_{24})$ and $\mathbf{b}=(b_{1},\ldots,b_{g})$
where $a_{j}\geq 0$, $b_{i}\geq 1$, and  $\sum a_{j}+\sum{b_{i}}=n+g$. The
corresponding divisor in $|C|$ is  
$$S+\sum_{i=1}^{g}b_{i}F_{i}+\sum_{j=1}^{24}a_{j}N_{j} $$
where $F_{i}$ is the smooth fiber containing $x_{i}$ and
$N_{1},\ldots,N_{24}$ are the nodal fibers.
\end{prop}

The proposition implies that $\mathcal{M}_{C,g}$ is the disjoint union of
components $\mathcal{M}_{\mathbf{a},\mathbf{b}}$ labeled by the vectors
$\mathbf{a}$ and $\mathbf{b}$.  
In section \ref{sec: analysis and local computations} we analyze the moduli
spaces $\mathcal{M}_{\mathbf{a},\mathbf{b}}$ in detail. The main result of
that section (Theorem \ref{thm: contribution of Mab}) is that the
contribution to $N_{g}(n)$ from $ 
\mathcal{M}_{\mathbf{a},\mathbf{b}}$ is the product of the local
contributions:
$$\prod_{i=1}^{g}b_{i}\sigma (b_{i})\prod_{j=1}^{24}p\left(
a_{j}\right) .$$ 

Our main theorem follows from this and some manipulations with the
generating functions. Recall that the
generating function of the 
partition function $p(l)$ is $\prod_{m=1}^{\infty }(1-q^{m})^{-1}$. 
Recall also that $\mathbf{a}$ is a $24$-tuple of integers with $a_{j}\geq 0$
and $\mathbf{b}$ is a $g$-tuple of integers with $b_{i}\geq 1$ and
$|\mathbf{a}|+|\mathbf{b}|=\sum_{j=1}^{24}a_{j}+\sum_{i=1}^{g} b_{i}=n+g$.
We compute:

\begin{eqnarray*}
N_{g}(n)q^{n}&=&	\left(\sum_{\begin{smallmatrix}
\mathbf{a},\mathbf{b}\\
|\mathbf{a}|+|\mathbf{b}|=n+g\end{smallmatrix}}\prod_{i=1}^{g}b_{i}\sigma
(b_{i})\prod_{j=1}^{24}p(a_{j}) \right)q^{n}\\
&=&\sum_{k=0}^{n}\sum_{\begin{smallmatrix}
|\mathbf{a}|=n-k\\
|\mathbf{b}|=g+k
\end{smallmatrix}}\prod_{i=1}^{g}b_{i}\sigma
(b_{i})q^{b_{i}-1}\prod_{j=1}^{24}p(a_{j}) q^{a_{j}}.
\end{eqnarray*}
Summing over $n$:
\begin{eqnarray*}
\sum_{n=0}^{\infty }N_{g}(n)q^{n}&=&\sum_{n=0}^{\infty
}\sum_{k=0}^{n}\left(\sum_{|\mathbf{b}|=g+k}\prod_{i=1}^{g}b_{i}\sigma
(b_{i})q^{b_{i}-1}
\right)\left(\sum_{|\mathbf{a}|=n-k}\prod_{j=1}^{24}p(a_{j})q^{a_{j}}
\right)\\  
&=&\sum_{k=0}^{\infty }\sum_{n=k}^{\infty
}\left(\sum_{|\mathbf{b}|=g+k}\prod_{i=1}^{g}b_{i}\sigma 
(b_{i})q^{b_{i}-1}
\right)\left(\sum_{|\mathbf{a}|=n-k}\prod_{j=1}^{24}p(a_{j})q^{a_{j}} \right)\\ 
&=&\left(\sum_{|\mathbf{b}|\geq g} \prod_{i=1}^{g}b_{i}\sigma
(b_{i})q^{b_{i}-1}\right)\left(\sum_{|\mathbf{a}|\geq
0}\prod_{j=1}^{24}p(a_{j})q^{a_{j}} \right)\\ 
&=&\prod_{i=1}^{g}\left(\sum_{b_{i}=1}^{\infty }b_{i}\sigma
(b_{i})q^{b_{i}-1} \right) \prod_{j=1}^{24}\left(\sum_{a_{j}=0}^{\infty
}p(a_{j}) q^{a_{j}} \right)\\
&=&\left(\sum_{b=1}^{\infty }b\sigma (b)q^{b-1}
\right)^{g}\prod_{m=1}^{\infty }(1-q^{m})^{-24}\\
&=&\left(\frac{d}{dq}G_{2}(q) \right)^{g}\frac{q}{\Delta (q)}.
\end{eqnarray*}

This proves our main theorem.

\section{Analysis of moduli spaces and local contributions}
\label{sec: analysis and local computations}      

The main goal of this section is to compute the contribution of the
component $\mathcal{M}_{\mathbf{a},\mathbf{b}}$ to the invariant
$N_{g}(n)$. Our strategy is simple in essence. We show that the moduli space
can be written as a product of various other moduli spaces and that the
tan\-gent-ob\-struct\-ion complex splits into factors that pull back from
tan\-gent-ob\-struct\-ion complexes on the other moduli spaces. We then show that
those individual moduli-obstruction problems have many components, each of
which can be identified with
moduli-obstruction problems arising for the Gro\-mov-Wit\-ten invariants of $\P
^{2}$ blown up multiple times. These contributions can then be determined
by elementary properties of the Gro\-mov-Wit\-ten invariants on blow ups of $\P
^{2}$. Using Cremona transformations, these contributions can be shown to all
either vanish or be equivalent to the number of straight lines between two
points (one). The computation then follows from straight forward
combinatorics.

This section is somewhat notationally heavy so to help the reader navigate we
summarize the notation used. We use
$\M _{C}$ for the full moduli space of stable genus $g$ maps to $X$ in the
class of $C$; $\M _{C,g}$ denotes the subset of $\M _{C}$ where the curves
pass through the $g$ points $x_{1},\ldots,x_{g}$. $\M _{C,g}$ breaks into
components $\M _{\mathbf{a},\mathbf{b}}$ indexed by vectors
$\mathbf{a}=(a_{1},\ldots,a_{24})$ and $\mathbf{b}=(b_{1},\ldots,b_{g})$
determining the image of the map. These components break into further
components indexed by ``data'' $\Lambda (\mathbf{b})$, $\lambda
(\mathbf{b})$, and $s(\mathbf{a})$ (Theorem \ref{thm: components of Mab}).
To prove Theorem \ref{thm: components of Mab}, it is convenient to
introduce $\M_{a}$ which denotes the moduli space of genus 0 stable maps to
$X$ with image $S+aN$ where $N$ is any fixed nodal fiber. The moduli space
$\M _{a}$ breaks into components because of the possibility of ``jumping''
behavior at the node of $N$. This behavior is encoded by certain kinds of
sequences 
$\{s_{n} \}$ (we call {\em admissible}) and hence the components of $\M
_{a}$ are indexed by such sequences. We denote those components by $\M
_{\{s_{n} \}}$. We compute the contribution of $\M _{\{s_{n} \}}$ by
``matching'' its tangent-obstruction complex with a tangent-obstruction
complex  on a moduli space of
stable maps to a blow up of $\P ^{2}$. This blow up is denoted $\til{P}$
and the relevant moduli space is denoted $\M ^{\til{P}}_{\{s_{n}
\}}$. Ultimately, we show that the contribution of each component 
to the invariant is either 0 or 1; those components that contribute 1 are
those for which the relavant admissable sequences have a special property
(we call such sequences {\em 1-admissable}). The contribution of $\M
_{\mathbf{a},\mathbf{b}}$ is then obtained by counting how many
possibilites there are for the data $\Lambda (\mathbf{b})$, $\lambda
(\mathbf{b})$, and $s(\mathbf{a})$ that have only 1-admissable sequences.

\subsection{Components of $\M _{\mathbf{a},\mathbf{b}} $}
We begin by identifying the connected components of
$\mathcal{M}_{\mathbf{a},\mathbf{b}}$. Call a sequence $\{s_{n} \}$ {\em
admissible} if each $s_{n}$ is a positive integer and the index $n$ runs
from some non-positive integer through some non-negative integer (the
sequence could consist solely of $\{s_{0} \}$ for example). Write $|s|$ for
$\sum_{n}s_{n}$. 

\begin{thm}\label{thm: components of Mab}
The connected components of $\mathcal{M}_{\mathbf{a},\mathbf{b}}$ are
indexed by the following data. For each $a_{i}\in \mathbf{a}$ assign an
admissible sequence $\{s_{n}(a_{i}) \}$ such that $|s(a_{i})|=a_{i}$. For each
$b_{j}\in \mathbf{b}$ assign a sublattice $\Lambda (b_{j})\subset \znums
\oplus  \znums $ of index $b_{j}$ and an element $\lambda (b_{j})$ of the
set $\{1,2,\ldots,b_{j} \}$. 
\end{thm}

\begin{rem}\label{rem: number of lattices is sigma}
The number of sublattices of $\znums \oplus \znums $ of index $b$ is
classically known and is given by $\sigma (b)=\sum_{d|b}d$. Thus we see
that the number of possible choices of data assigned to $\mathbf{b}$ is
$\prod_{i=1}^{g}b_{i}\sigma (b_{i})$. 
\end{rem}

Let $f\in \M _{\mathbf{a},\mathbf{b}}$. Since the image of $f:D\to X$ is
$S+\sum_{i}b_{i}F_{i}+\sum_{j}a_{j}N_{j}$ and is reducible, $D$ must be
reducible and its components must group into the set of components mapping
to $S$, $F_{1},\ldots,F_{g}$, and $N_{1},\ldots,N_{24}$. Since the
components mapping to each $F_{i}$ must have geometric genus at least 1 and
the total geometric genus of $D$ is $g$, the components of $D$ mapping
to $F_{i}$ must each be genus 1 and all other components of $D$ are
rational. Furthermore, the dual graph of $D$ is a tree and the $g$ marked
points are on the $g$ elliptic components which we call
$G_{1},\ldots,G_{g}$. Then since $f$ is a stable map, the image of all
ghost components 
of $D$ must lie in the nodal fibers $N_{1},\ldots,N_{24}$. We denote the
component of $D$ mapping isomorphically onto $S$ also by $S$. 

So far then, we can describe the domain $D$ as a rational curve $S$ that
has attached to it $g$ marked elliptic curves $G_{1},\ldots,G_{g}$ and 24
components $D_{1},\ldots,D_{24}$ that are either empty (if $a_{i}=0$) or a
tree of rational components. Furthermore, $f|_{G_{i}}:G_{i}\to F_{i}$ is a
degree $b_{i}$ map preserving the intersection with $S$ and sending the
marked point to $x_{i}$ and $f|D_{j}:D_{j}\to N_{j}$ has total degree
$a_{j}$. 

Using the intersection with $S$ as an origin for $G_{i}$ and $F_{i}$, we
can identify the number of distinct possibilities for the map $f:G_{i}\to
F_{i}$ with the number of degree $b_{i}$ homomorphisms onto a fixed
elliptic curve $F_{i}$. This is precisely the number of index $b_{i}$
sublattices of $\znums \oplus \znums $. Additionally, since $x_{i}$ has
$b_{i}$ preimages under $f$ (the $x_{i}$'s are chosen generically), there are
$b_{i}$ choices for the location of the marked point on $G_{i}$ for each
homomorphism  $f:G_{i}\to F_{i}$. Thus the
data associated to $\mathbf{b}$ in the theorem completely determines $f$
restricted to $G_{1},\ldots,G_{g}$. Now $f|_{S}$ is determined and so we
can reconstruct $f$ completely from $f|_{N_{1}},\ldots,f|_{N_{24}}$ and the
data $\Lambda (b_{i})$, $\lambda (b_{i})$. It follows that the subset of
$\M _{\mathbf{a},\mathbf{b}}$ with fixed data for $\mathbf{b}$ is the
product of the moduli spaces $\prod_{j=1}^{24}\M _{[a_{j}]_{j},0}$, where
$[c]_{j}$ denotes the 24-tuple $(0,\ldots,c,\ldots,0)$ with $c$ in
the $i$th slot and zeros elsewhere. The connected components of $\M
_{\mathbf{a},\mathbf{b}}$ are in one to one correspondence with the data
$\Lambda (b_{i})$, $\lambda (b_{i})$, and $s(a_{j})$ and Theorem \ref{thm:
components of Mab} is proved provided we can show that the connected
components of $\M _{[a],0}$ are in one to one correspondence with
admissible sequences $s$ of magnitude $|s|=a$. We state this as a   
\begin{lem}\label{lem: components of nodal multiple covers are admissable
seqs} 
Let $\M _{a}$ be the moduli space of stable, genus 0 maps to $X$ with image
$S+aN$ for any fixed nodal fiber $N$. Then $\M _{a}$ is a disjoint union
$\coprod _{\{s_{n} \}}\M 
_{\{s_{n} \}}$ of spaces $\M _{\{s_{n} \}}$ labeled by admissible
sequences $\{s_{n} \}$ with $|s|=\sum_{n}s_{n}=a$. 
\end{lem}
{\sc Proof: }
Let $\Sigma(a)$ be a genus 0 nodal curve consisting of a linear chain of
$2a+1$ smooth components $\Sigma _{-a},\ldots,\Sigma _{a}$ with an
additional 
component $\Sigma _{*}$ meeting $\Sigma _{0}$ (so $\Sigma_{n} \cap\Sigma
_{m}= \emptyset   $ unless $|n-m|= 1$ and $\Sigma _{*}\cap \Sigma
_{n}=\emptyset $ unless $n=0$). Fix a map of $\Sigma(a) $ to $X$ with image
$S\cup N$ in the following way. Map $\Sigma _{*}$ to $S$ with degree 1
and map each $\Sigma _{n}$ to $N$ with degree one. Require that a
neighborhood of each singular point $\Sigma _{n}\cap \Sigma _{n+1}$ is
mapped biholomorphically onto its image with $\Sigma _{n}\cap \Sigma
_{n+1}$ mapping to  the nodal point of $N$. 

Let $\{s_{n} \}$ be an admissible sequence with $|s|=a$. Since the index
$n$ of the sequence cannot be smaller than $-a$ or larger than $a$, we can
extend $\{s_{n} \}$ to a sequence $s_{-a},\ldots,s_{a}$ by setting
$s_{n}=0$ for those not previously defined. Define $\M _{\{s_{n} \}}$ to be
the moduli space of genus 0 stable maps {\em to} $\Sigma(a) $ in the class
$$
\Sigma _{*} +\sum_{n=-a}^{a}s_{n}\Sigma _{n}.
$$
By composition with the fixed map from $\Sigma(a) $ to $X$, we get a stable
map in $\M _{a}$ from each map in $\M _{\{s_{n} \}}$. To prove the lemma we
need to show\footnote{To prove the
lemma as stated we also need to
show that $\M _{\{s_{n} \}}$ is connected. This is not hard, but
since we never actually use this part of the result, we will leave its
proof to the reader.} that 
every map in $\M _{a}$ factors uniquely in this way through a map in $\M
_{\{s_{n} \}}$ for some admissible $\{s_{n} \}$. 
The following figure illustrates some of the phenomenon that can occur. The
numbers on components of $D$ indicate the degree of $f$ on that
component. The $A$'s and $B$'s indicate the local behavior of the map when
a nodal point of $D$ is mapped to the nodal point in $N$. 
 
\vskip 15pt
\setlength{\unitlength}{0.00083333in}
\begingroup\makeatletter\ifx\SetFigFont\undefined
% extract first six characters in \fmtname
\def\x#1#2#3#4#5#6#7\relax{\def\x{#1#2#3#4#5#6}}%
\expandafter\x\fmtname xxxxxx\relax \def\y{splain}%
\ifx\x\y   % LaTeX or SliTeX?
\gdef\SetFigFont#1#2#3{%
  \ifnum #1<17\tiny\else \ifnum #1<20\small\else
  \ifnum #1<24\normalsize\else \ifnum #1<29\large\else
  \ifnum #1<34\Large\else \ifnum #1<41\LARGE\else
     \huge\fi\fi\fi\fi\fi\fi
  \csname #3\endcsname}%
\else
\gdef\SetFigFont#1#2#3{\begingroup
  \count@#1\relax \ifnum 25<\count@\count@25\fi
  \def\x{\endgroup\@setsize\SetFigFont{#2pt}}%
  \expandafter\x
    \csname \romannumeral\the\count@ pt\expandafter\endcsname
    \csname @\romannumeral\the\count@ pt\endcsname
  \csname #3\endcsname}%
\fi
\fi\endgroup
{\renewcommand{\dashlinestretch}{30}
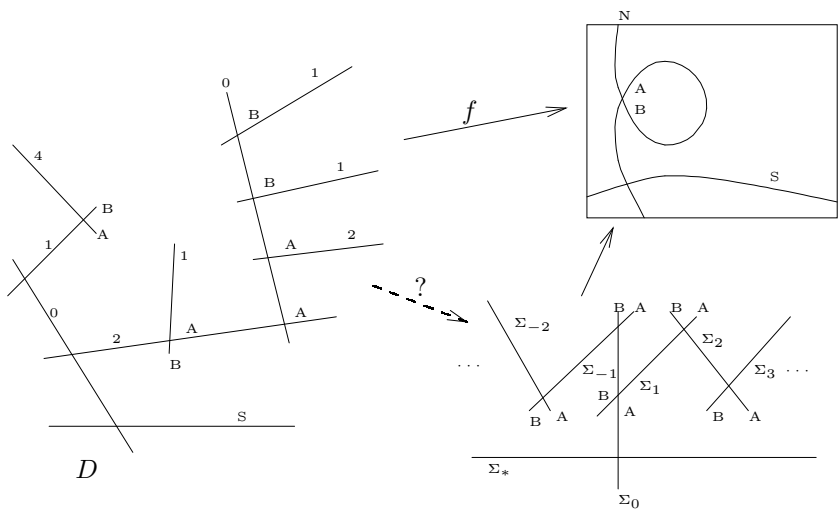
\begin{picture}(5199,3100)(0,-10)
\path(272,499)(1802,499)
\path(793,336)(45,1541)
\path(12,1313)(565,1867)
\path(565,1703)(45,2256)
\path(240,922)(2062,1183)
\path(1020,954)(1054,1638)
\path(1769,1020)(1379,2582)
\path(1345,2256)(2160,2745)
\path(1444,1900)(2323,2094)
\path(1543,1541)(2355,1638)
\path(3625,3005)(5187,3005)(5187,1801)
	(3625,1801)(3625,3005)
\path(2908,304)(5057,304)
\path(3820,109)(3820,1214)
\path(3690,564)(4309,1183)
\path(4145,1214)(4634,596)
\path(4373,596)(4894,1183)
\path(3918,1214)(3267,596)
\path(3592,1313)(3788,1737)
\path(3764.879,1615.487)(3788.000,1737.000)(3710.417,1640.663)
\path(2486,2289)(3494,2485)
\path(3381.932,2432.647)(3494.000,2485.000)(3370.480,2491.544)
\dashline{60.000}(2291,1377)(2876,1150)
\path(2753.275,1165.442)(2876.000,1150.000)(2774.980,1221.379)
\path(3397,596)(3006,1281)
\path(3625,1931)	(3674.856,1947.727)
	(3721.298,1963.087)
	(3764.506,1977.129)
	(3804.656,1989.903)
	(3876.500,2011.839)
	(3938.259,2029.286)
	(3991.360,2042.635)
	(4037.232,2052.276)
	(4077.303,2058.601)
	(4113.000,2062.000)

\path(4113,2062)	(4150.433,2063.312)
	(4189.990,2063.372)
	(4232.064,2062.122)
	(4277.048,2059.502)
	(4325.335,2055.456)
	(4377.320,2049.924)
	(4433.394,2042.847)
	(4493.951,2034.167)
	(4559.385,2023.826)
	(4630.089,2011.765)
	(4667.540,2005.071)
	(4706.456,1997.926)
	(4746.886,1990.320)
	(4788.879,1982.249)
	(4832.484,1973.703)
	(4877.752,1964.676)
	(4924.729,1955.161)
	(4973.467,1945.150)
	(5024.014,1934.635)
	(5076.419,1923.610)
	(5130.731,1912.068)
	(5187.000,1900.000)

\path(3984,1801)	(3946.478,1870.787)
	(3914.802,1931.335)
	(3888.481,1983.705)
	(3867.024,2028.957)
	(3836.733,2102.351)
	(3820.000,2160.000)

\path(3820,2160)	(3808.301,2231.199)
	(3803.880,2274.550)
	(3801.209,2320.276)
	(3800.886,2366.280)
	(3803.512,2410.466)
	(3820.000,2485.000)

\path(3820,2485)	(3838.459,2527.277)
	(3861.845,2574.734)
	(3890.291,2623.974)
	(3923.926,2671.599)
	(3962.884,2714.212)
	(4007.294,2748.415)
	(4057.289,2770.810)
	(4113.000,2778.000)

\path(4113,2778)	(4166.998,2770.225)
	(4217.344,2751.438)
	(4262.791,2723.024)
	(4302.096,2686.366)
	(4334.014,2642.851)
	(4357.301,2593.861)
	(4370.711,2540.783)
	(4373.000,2485.000)

\path(4373,2485)	(4364.917,2436.220)
	(4347.505,2391.179)
	(4289.399,2316.354)
	(4251.058,2288.589)
	(4208.093,2268.602)
	(4161.682,2257.403)
	(4113.000,2256.000)

\path(4113,2256)	(4051.436,2271.210)
	(3998.505,2302.663)
	(3953.453,2346.476)
	(3915.526,2398.765)
	(3883.971,2455.646)
	(3858.034,2513.234)
	(3836.962,2567.647)
	(3820.000,2615.000)

\path(3820,2615)	(3805.617,2675.026)
	(3802.021,2711.564)
	(3800.823,2753.967)
	(3802.021,2803.375)
	(3805.617,2860.925)
	(3811.610,2927.754)
	(3815.505,2965.004)
	(3820.000,3005.000)

\put(2844,2419){\makebox(0,0)[lb]{\smash{{{\SetFigFont{10}{12.0}{rm}$f$}}}}}
\put(3820,3038){\makebox(0,0)[lb]{\smash{{{\SetFigFont{5}{6.0}{rm}N}}}}}
\put(4764,2030){\makebox(0,0)[lb]{\smash{{{\SetFigFont{5}{6.0}{rm}S}}}}}
\put(3918,2582){\makebox(0,0)[lb]{\smash{{{\SetFigFont{5}{6.0}{rm}A}}}}}
\put(3918,2452){\makebox(0,0)[lb]{\smash{{{\SetFigFont{5}{6.0}{rm}B}}}}}
\put(3690,662){\makebox(0,0)[lb]{\smash{{{\SetFigFont{5}{6.0}{rm}B}}}}}
\put(3854,564){\makebox(0,0)[lb]{\smash{{{\SetFigFont{5}{6.0}{rm}A}}}}}
\put(3788,1214){\makebox(0,0)[lb]{\smash{{{\SetFigFont{5}{6.0}{rm}B}}}}}
\put(3918,1214){\makebox(0,0)[lb]{\smash{{{\SetFigFont{5}{6.0}{rm}A}}}}}
\put(4145,1214){\makebox(0,0)[lb]{\smash{{{\SetFigFont{5}{6.0}{rm}B}}}}}
\put(4309,1214){\makebox(0,0)[lb]{\smash{{{\SetFigFont{5}{6.0}{rm}A}}}}}
\put(4634,532){\makebox(0,0)[lb]{\smash{{{\SetFigFont{5}{6.0}{rm}A}}}}}
\put(3006,206){\makebox(0,0)[lb]{\smash{{{\SetFigFont{5}{6.0}{rm}$\Sigma_*$}}}}}
\put(3820,11){\makebox(0,0)[lb]{\smash{{{\SetFigFont{5}{6.0}{rm}$\Sigma_0$}}}}}
\put(3950,727){\makebox(0,0)[lb]{\smash{{{\SetFigFont{5}{6.0}{rm}$\Sigma_1$}}}}}
\put(4340,1020){\makebox(0,0)[lb]{\smash{{{\SetFigFont{5}{6.0}{rm}$\Sigma_2$}}}}}
\put(4666,825){\makebox(0,0)[lb]{\smash{{{\SetFigFont{5}{6.0}{rm}$\Sigma_3$}}}}}
\put(4406,532){\makebox(0,0)[lb]{\smash{{{\SetFigFont{5}{6.0}{rm}B}}}}}
\put(3267,499){\makebox(0,0)[lb]{\smash{{{\SetFigFont{5}{6.0}{rm}B}}}}}
\put(3430,532){\makebox(0,0)[lb]{\smash{{{\SetFigFont{5}{6.0}{rm}A}}}}}
\put(4862,825){\makebox(0,0)[lb]{\smash{{{\SetFigFont{5}{6.0}{rm}$\cdots$}}}}}
\put(3592,825){\makebox(0,0)[lb]{\smash{{{\SetFigFont{5}{6.0}{rm}$\Sigma_{-1}$}}}}}
\put(2551,1313){\makebox(0,0)[lb]{\smash{{{\SetFigFont{10}{12.0}{rm}?}}}}}
\put(1444,532){\makebox(0,0)[lb]{\smash{{{\SetFigFont{5}{6.0}{rm}S}}}}}
\put(1020,858){\makebox(0,0)[lb]{\smash{{{\SetFigFont{5}{6.0}{rm}B}}}}}
\put(1119,1084){\makebox(0,0)[lb]{\smash{{{\SetFigFont{5}{6.0}{rm}A}}}}}
\put(1802,1183){\makebox(0,0)[lb]{\smash{{{\SetFigFont{5}{6.0}{rm}A}}}}}
\put(1738,1606){\makebox(0,0)[lb]{\smash{{{\SetFigFont{5}{6.0}{rm}A}}}}}
\put(1607,1996){\makebox(0,0)[lb]{\smash{{{\SetFigFont{5}{6.0}{rm}B}}}}}
\put(1509,2419){\makebox(0,0)[lb]{\smash{{{\SetFigFont{5}{6.0}{rm}B}}}}}
\put(1345,2615){\makebox(0,0)[lb]{\smash{{{\SetFigFont{5}{6.0}{rm}0}}}}}
\put(1900,2680){\makebox(0,0)[lb]{\smash{{{\SetFigFont{5}{6.0}{rm}1}}}}}
\put(2062,2094){\makebox(0,0)[lb]{\smash{{{\SetFigFont{5}{6.0}{rm}1}}}}}
\put(2128,1671){\makebox(0,0)[lb]{\smash{{{\SetFigFont{5}{6.0}{rm}2}}}}}
\put(663,1020){\makebox(0,0)[lb]{\smash{{{\SetFigFont{5}{6.0}{rm}2}}}}}
\put(1086,1541){\makebox(0,0)[lb]{\smash{{{\SetFigFont{5}{6.0}{rm}1}}}}}
\put(597,1833){\makebox(0,0)[lb]{\smash{{{\SetFigFont{5}{6.0}{rm}B}}}}}
\put(565,1671){\makebox(0,0)[lb]{\smash{{{\SetFigFont{5}{6.0}{rm}A}}}}}
\put(175,2160){\makebox(0,0)[lb]{\smash{{{\SetFigFont{5}{6.0}{rm}4}}}}}
\put(240,1606){\makebox(0,0)[lb]{\smash{{{\SetFigFont{5}{6.0}{rm}1}}}}}
\put(272,1183){\makebox(0,0)[lb]{\smash{{{\SetFigFont{5}{6.0}{rm}0}}}}}
\put(3169,1118){\makebox(0,0)[lb]{\smash{{{\SetFigFont{5}{6.0}{rm}$\Sigma_{-2}$}}}}}
\put(2811,858){\makebox(0,0)[lb]{\smash{{{\SetFigFont{5}{6.0}{rm}$\cdots$}}}}}
\put(435,174){\makebox(0,0)[lb]{\smash{{{\SetFigFont{10}{12.0}{rm}$D$}}}}}
\end{picture}
}
\vskip 15pt

Consider the dual graph of the domain $D$ of a map $f:D\to X$ in $\M
_{a}$. The graph is a tree with one special vertex $v_{*}$ (the component
mapping to $S$) whose valence is 1. Every other vertex $v$ is marked with a
non-negative integer $l_{v}$ (the degree of the component associated to
$v$) such that the sum of the $l_{v}$'s is $a$. Vertices with a marking of
0 (ghost components) must have valence at least three (stability). We mark
the edges in the following way. Each edge corresponds to a nodal
singularity in $D$ and if the node is not mapped to the nodal point in $N$,
we do not mark the edge. The remaining edges are marked with either a pair
of the letters $A$ or $B$, a single letter of $A$ or $B$, or nothing as
follows. Label the two branches near the node in $N$ by $A$ and $B$. Three
things can then happen for an edge corresponding to a node in $D$ that gets
mapped to the node in $N$.\begin{enumerate}
\item If the edge connects two ghost components, do not mark the edge.
\item If the edge connects one ghost component with one non-ghost
component, mark the edge with an $A$ or $B$ depending on whether the
non-ghost component is mapped (locally) to the $A$ branch or the $B$ branch.
\item Finally, if the edge connects two non-ghost components, then mark the
edge with two of the letters $A$ or $B$, one near each of the vertices,
according to which branch that corresponding component maps to
(locally). Note that all the combinations $AB$, $BA$, $AA$, and $BB$ can
occur. 
\end{enumerate}
The markings on our graph now tell us how and when the map ``jumps''
branches. Jumping from $A$ to $B$ will correspond to moving from $\Sigma
_{n}$ to $\Sigma _{n+1}$ in the factored map. To determine which component
$\Sigma _{n}$ a component of $D$ gets mapped to, we count how many
``jumps'' occur between the corresponding vertex $v$ and the central vertex
$v_{*}$. For every non-ghost component vertex $v$ assign its index $n_{v}$
by traveling from $v_{*}$ to $v$ in the graph and counting $+1$ for each
$AB$ pair passed through, $-1$ for each $BA$ pair, and $0$ for each $AA$ or
$BB$ pair. We can now uniquely factor $f:D\to X$ through the fixed map
$\Sigma (a)\to X$. The component of $D$ corresponding to a vertex $v$ gets
mapped to $\Sigma _{n_{v}}$. The factorization is unique since away from
the $AB$ or $BA$ jumps, $f$ factors uniquely through the normalization of
$N$. \qed 

The marked dual graph for the previously illustrated example is
below.
\vskip 15pt
\setlength{\unitlength}{0.00083333in}
\begingroup\makeatletter\ifx\SetFigFont\undefined
% extract first six characters in \fmtname
\def\x#1#2#3#4#5#6#7\relax{\def\x{#1#2#3#4#5#6}}%
\expandafter\x\fmtname xxxxxx\relax \def\y{splain}%
\ifx\x\y   % LaTeX or SliTeX?
\gdef\SetFigFont#1#2#3{%
  \ifnum #1<17\tiny\else \ifnum #1<20\small\else
  \ifnum #1<24\normalsize\else \ifnum #1<29\large\else
  \ifnum #1<34\Large\else \ifnum #1<41\LARGE\else
     \huge\fi\fi\fi\fi\fi\fi
  \csname #3\endcsname}%
\else
\gdef\SetFigFont#1#2#3{\begingroup
  \count@#1\relax \ifnum 25<\count@\count@25\fi
  \def\x{\endgroup\@setsize\SetFigFont{#2pt}}%
  \expandafter\x
    \csname \romannumeral\the\count@ pt\expandafter\endcsname
    \csname @\romannumeral\the\count@ pt\endcsname
  \csname #3\endcsname}%
\fi
\fi\endgroup
{\renewcommand{\dashlinestretch}{30}
\begin{picture}(4789,2063)(0,-10)
\path(48,986)(873,647)(1600,647)
	(2570,986)(2327,1809)
\path(1600,647)(1648,65)
\path(2570,986)(3394,986)
\path(3394,986)(4364,1374)
\path(3443,1858)(3394,986)(3539,211)
\put(0,1083){\makebox(0,0)[lb]{\smash{{{\SetFigFont{7}{8.4}{rm}(4,-1)}}}}}
\put(97,792){\makebox(0,0)[lb]{\smash{{{\SetFigFont{7}{8.4}{rm}A}}}}}
\put(582,647){\makebox(0,0)[lb]{\smash{{{\SetFigFont{7}{8.4}{rm}B}}}}}
\put(825,744){\makebox(0,0)[lb]{\smash{{{\SetFigFont{7}{8.4}{rm}(1,0)}}}}}
\put(1552,744){\makebox(0,0)[lb]{\smash{{{\SetFigFont{7}{8.4}{rm}0}}}}}
\put(1745,17){\makebox(0,0)[lb]{\smash{{{\SetFigFont{9}{8.4}{rm}$v_*$}}}}}
\put(2521,840){\makebox(0,0)[lb]{\smash{{{\SetFigFont{7}{8.4}{rm}(2,0)}}}}}
\put(2909,1034){\makebox(0,0)[lb]{\smash{{{\SetFigFont{7}{8.4}{rm}A}}}}}
\put(2521,1180){\makebox(0,0)[lb]{\smash{{{\SetFigFont{7}{8.4}{rm}A}}}}}
\put(2425,1568){\makebox(0,0)[lb]{\smash{{{\SetFigFont{7}{8.4}{rm}B}}}}}
\put(2279,1858){\makebox(0,0)[lb]{\smash{{{\SetFigFont{7}{8.4}{rm}(1,1)}}}}}
\put(3345,1955){\makebox(0,0)[lb]{\smash{{{\SetFigFont{7}{8.4}{rm}(2,0)}}}}}
\put(4461,1374){\makebox(0,0)[lb]{\smash{{{\SetFigFont{7}{8.4}{rm}(1,1)}}}}}
\put(3636,162){\makebox(0,0)[lb]{\smash{{{\SetFigFont{7}{8.4}{rm}(1,1)}}}}}
\put(3491,889){\makebox(0,0)[lb]{\smash{{{\SetFigFont{7}{8.4}{rm}0}}}}}
\put(3491,1568){\makebox(0,0)[lb]{\smash{{{\SetFigFont{7}{8.4}{rm}A}}}}}
\put(4025,1083){\makebox(0,0)[lb]{\smash{{{\SetFigFont{7}{8.4}{rm}B}}}}}
\put(3539,453){\makebox(0,0)[lb]{\smash{{{\SetFigFont{7}{8.4}{rm}B}}}}}
\end{picture}
}
\vskip 15pt
Here we've marked the vertices with $(l_{v},n_{v})$ so in this example
$s_{-1}=4$, $s_{0}=5$, and $s_{1}=3$.

\subsection{The tangent-obstruction complex}

To compute the contribution of $\M _{\mathbf{a},\mathbf{b}}$ to
$N_{g}(n)=\Psi (C,g,g)$, we recall the definition of the invariant. Let $\M
_{C}$ be the moduli space of $g$-marked, genus $g$ stable maps to $X$ (with
its twistor family) in the class of $C$. Let $\M _{C,g}\subset \M _{C}$
denote the restriction to those maps that send the $i$th marked point to the
point $x_{i}\in X$; the virtual dimension of $\M _{C,g}$ is 0. By
definition, $N_{g}(n) $ is the evaluation of $ev^{*}_{1}([x_{1}])\wedge
\cdots \wedge  ev^{*}_{g}([x_{g}])$ on $[\M _{C}]^{vir}$ which is the same
as the class $[\M _{C,g}]^{vir}$. 

Li and Tian construct the class $[\M _{C,g}]^{vir}$ in the symplectic
category in \cite{Li-Tian} and in the algebraic category in
\cite{Li-Tian2}. In \cite{Li-Tian3} they show the classes coincide when the
symplectic manifold is algebraic. This enables us to compute 
purely algebro-geometrically although we need the machinery of the
symplectic category to define the invariant.

The data of the algebraic construction of $[\M _{C,g}]^{vir}$ is the
tan\-gent-ob\-struct\-ion complex $[T^{1}\to T^{2}]$. It is a complex of sheaves
over $\M _{C,g}$ whose stalks over a map $f:D\to X$ fit into an exact
sequence (c.f. \cite{Gr-Pa}) 
$$
\begin{CD}
0@>>>\operatorname{Aut}(D)@>>>H^{0}(D,f^{*}(TX))@>>>T^{1}@>>>\\
@>>>\operatorname{Def}(D)\oplus T_{t_{0}
}B@>>>H^{1}(D,f^{*}(TX))@>>>T^{2}@>>>0 
\end{CD}
$$
where $\operatorname{Aut}(D)$ is the space of infinitesimal automorphisms
of the domain, $\operatorname{Def}(D)$ is the space of infinitesimal
deformations of the domain, and $T_{t_{0} }B$ is the tangent space of the
family of K\"ahler structures at $t_{0}$ (formally the pull back by
$J_{B}:B\to \mathcal{J} $ of the tangent space of $\mathcal{J}$ restricted
to $t_{0}\in B$). The sequence can be informally interpreted as
follows. The tangent 
space to $\M _{C,g}$ at a map $f:D\to X$ contains those vector fields along
the map modulo those that come from vector fields on the domain. The
tangent space also contains infinitesimal deformations of the complex
structure on the 
domain and infinitesimal deformations of the complex structure of the range
in the family, however, those infinitesimal deformations are obstructed
from being actual deformations if they have non-zero image in
$H^{1}(D,f^{*}(TX))$. 

If the ranks of $T^{i}$ remain constant over the moduli space, then the
moduli space is a smooth orbifold and $T^{2}$ is a smooth orbi-bundle. In
this case, the virtual fundamental class is simply the Euler class
of $T^{2}$. This is, in fact, the case for $\M _{C,g}$. Since we have
identified the components of $\M _{C,g}$ with various products of the
moduli spaces $\M _{\{s_{n} \}}$ which are smooth orbifolds ($\M _{\{s_{n}
\}}$ is smooth since it is a moduli space of stable maps to a range
consisting of convex varieties meeting at general points). In the rest of
this section we show how the bundle $T^{2}$ can be expressed as a sum of
bundles pulling back from the factors of the product.

The exact sequence shows that when (a component of) the moduli space of
stable 
maps all have the same image, the tangent obstruction complex only depends
upon the restriction of the tangent bundle to $X$ to the image of $f$ (and
not on the rest of $X$). 

Breaking our moduli space into connected components, we have 
$$[\M _{C,g}]^{vir}=\sum_{\mathbf{a},\mathbf{b}}\sum_{s,\Lambda
, \lambda} [\prod_{i=1}^{24}\M
_{s(a_{i})}]^{vir}.$$ 
We wish to determine $[\prod_{i=1}^{24}\M _{s(a_{i})}]^{vir}$ in terms of
tan\-gent-ob\-struct\-ion complexes on the $\M _{s(a_{i})}$'s.
To specify a tan\-gent-ob\-struct\-ion complex for $\M _{s(a_{i})}$ it is enough
to define a rank 2 bundle $T$ on $\Sigma(a_{i}) $ and use the above sequence to
define $T^{1}$ and $T^{2}$. Let $T$ be the bundle on $\Sigma(a_{i}) $ that
restricted to each component $\Sigma _{i}$ is $T\Sigma _{i}\oplus
\mathcal{O}_{\Sigma _{i}}(-2)$ and restricted to $\Sigma _{*}$ is $T\Sigma
_{*}\oplus \mathcal{O}_{\Sigma _{*}}(-1)$. Call the resulting
tan\-gent-ob\-struct\-ion complex $[T^{1}\to T^{2}]_{s(a_{i})}$. Note that the
bundle obtained by pulling back $TX$ by 
the fixed map $\Sigma(a_{i}) \to S\cup N_{i}$ is isomorphic to $T\otimes 
\mathcal{O}_{\Sigma _{*}}(1)$. 

In this subsection we show:
\begin{lem}\label{lem: t-o complex is product of t-o complexs from fibers}
The tan\-gent-ob\-struct\-ion complex defining $[\prod_{i=1}^{24}\M
_{s(a_{i})}]^{vir}$ is isomorphic to the direct sum of tan\-gent-ob\-struct\-ion
complexes $[T^{1}\to T^{2}]_{s(a_{i})}$.
\end{lem}

In the case at hand, $B$ is the hyperk\"ahler family
$S(\mathcal{H}^{2}_{+,g})$ and 
so the tangent space to $t_{0}$ is the space perpendicular to $\omega
_{t_{0}}$ in $\mathcal{H}^{2}_{+,g}$ which can be canonically identified
with the space of holomorphic 2-forms. Thus $T_{t_{0}}B$ is
canonically\footnote{There is a question how
to orient the twistor family which corresponds to whether we identify the
perpendicular to $\omega _{0}$ with $H^{0}(X,\mathcal{O})$ or
$H^{0}(X,\mathcal{O})^{*}$. We choose the convention that makes the
intersection of the twistor family with the set of projective $K3$'s
positive.} 
$H^{0}(X,K)^{*}\cong H^{2}(X,\mathcal{O})$. 
It is convienent to use the exact sequence
$$0\to \mathcal{O}\to \mathcal{O}(S)\to N_{S}\to 0 $$
to identify $H^{2}(X,\mathcal{O})$ with $H^{1}(X,N_{S})$ where $N_{S} $ is
the normal bundle to $S$ in $X$. Since $TS$ is positive we can in fact
identify $T_{t_{0}}B$ with $H^{1}(S,TX|_{S})$. 

We know that infinitesimal deformations of $f$ in the direction of the
twistor family are obstructed (Proposition \ref{prop: one curve in class of
C in twistor family}), and we wish to determine the image of $T_{t_{0}}B\to
H^{1}(f^{*}(TX))$. Let $D'$ be the union of the components
$D_{1},\ldots,D_{24}$ and $G_{1},\ldots,G_{g}$ in $D$ so that $D=D'\cup
S$. Twisting the (partial) normalization sequence
$$0\to \mathcal{O}_{D}\to \mathcal{O}_{S}\oplus \mathcal{O}_{D'}\to
\mathcal{O}_{S\cap D'}\to 0 $$ 
by $f^{*}(TX)$ and taking cohomology, we see that $H^{1}(D,f^{*}(TX))$
surjects onto $H^{1}(S,f^{*}(TX))$. The composition 
$$T_{t_{0}}B\cong H^{1}(S,f^{*}(TX))\to H^{1}(D,f^{*}(TX))\to
H^{1}(S,f^{*}(TX))  $$
is the identity. We can thus rewrite the sequence for $T^{1}$ and $T^{2}$
as 
\begin{eqnarray}\label{eqn: sequence for t-o complex without familiy term}
%\begin{CD}
%0@>>>\operatorname{Aut}(D)@>>>H^{0}(D,f^{*}(TX))@>>>T^{1}@>>>\\
%@>>>\operatorname{Def}(D)
%@>>>H^{1}(D,f^{*}(TX))/H^{1}(S,f^{*}(TX))@>>>T^{2}@>>>0  .
%\end{CD}
0\to \operatorname{Aut}(D)\to & H^{0}(D,f^{*}(TX))&\to T^{1}\to \\
\nonumber \to \operatorname{Def}(D)
\to & H^{1}(D,f^{*}(TX))/H^{1}(S,f^{*}(TX))&\to T^{2}\to 0  .
\end{eqnarray}
This identifies our tan\-gent-ob\-struct\-ion complex with an equivalent
tan\-gent-ob\-struct\-ion complex without a family: since the tan\-gent-ob\-struct\-ion
complex only depends upon the restriction of $TX$ to the image of $f$, we
can define a new tan\-gent-ob\-struct\-ion complex over $\M
_{\mathbf{a},\mathbf{b}}$ by dropping the dependence on the family and
twisting the normal bundle of $S$ by $\mathcal{O}(1)$ so that
$N_{S}=\mathcal{O}_{S}(-1)$. The above exact sequence shows that the
resulting complex is isomorphic to $[T^{1}\to T^{2}]$ since the
effect of changing the normal bundle of $S$ from $\mathcal{O}(-2)$ to
$\mathcal{O}(-1)$ is solely to replace $H^{1}(D,f^{*}(TX))$ by
$H^{1}(D,f^{*}(TX))/H^{1}(S,\mathcal{O}(-2))$. 

\begin{rem}\label{rem: obstr argumnet applies to E(1) also}
The discussion of this section also applies to the rational elliptic
surface where we consider ordinary Gro\-mov-Wit\-ten invariants (no family) but
the normal bundle to the section is already degree $-1$.
\end{rem}

To complete the proof of the lemma we now use the (partial) normalization
sequence 
$$0\to \mathcal{O}_{D}\to \mathcal{O}_{S}\oplus _{i}
\mathcal{O}_{G_{i}}\oplus _{j} \mathcal{O}_{D_{j}}\to  \oplus _{i}
\mathcal{O}_{S\cap G_{i}}\oplus _{j} \mathcal{O}_{S\cap D_{j}}\to 0 $$ 
twisted by $f^{*}(TX)$.
With care one can use the isomorphisms obtained from the cohomology
sequence to see that the Sequence \ref{eqn: sequence for t-o complex
without familiy term} is a direct sum of the sequences 
defining $[T^{1}\to T^{2}]_{s(a_{i})}$.
The dependence on the $G_{j}$ components goes away: the $TX|_{S\cap G_{i}}$
term arising from the normalization sequence cancels with the infinitesimal
deformation of $D$ that smooths the intersection of $S$ and $G_{i}$ and the
new automorphism of $S$ that occurs when the intersection with $G_{i}$ is
removed. Finally, we use the fixed maps $\Sigma(a_{i}) \to S\cup N_{i}$ to
identify the $H^{*}(D,f^{*}(TX))$ terms with those terms
$H^{*}(D,\til{f}^{*}(T))$ coming from maps $\til{f}:D\to \Sigma(a_{i}) $ in
$\M_{s(a_{i})}$.\qed 

\subsection{Computations via blow-ups on $\P ^{2}$}
In the previous subsection, we showed that the virtual fundamental cycle of
a component $\prod_{i=1}^{24}\M _{s(a_{i})}$ (for any fixed
$s(\mathbf{a})$, $\Lambda
(\mathbf{b})$, and $\lambda (\mathbf{b})$) is given by the product of
virtual 
fundamental cycles on $\M _{s(a_{i})}$ defined by the
tan\-gent-ob\-struct\-ion 
complex $[T^{1}\to T^{2}]_{s(a_{i})}$. In this subsection, we will realize
the 
moduli-obstruction problem $(\M _{s(a_{i})},[T^{1}\to T^{2}]_{s(a_{i})})$
as one coming from $\til{P}$, a certain blow-up of $\P ^{2}$ at $2a+3$
points.

 The homology classes of $\til{P}$ will have a diagonal basis $h$,
$e_{-a-1},\ldots,e_{a+1}$ where $h^{2}=1$ and $e_{n}^{2}=-1$.
We construct $\til{P}$ as follows. Begin with a linear $\cnums ^{*}$ action
on $\P ^{2}$ fixing a line $H$ and a point $p$. Choose three
points $p_{-}$, $p_{0}$, and $p_{+}$ on $H$ and blow them up to
obtain three exceptional curves $E_{-1}$, $E_{0}$, and $E_{1}$ representing
classes $e_{-1}$, $e_{0}$, and $e_{1}$. The proper transform of
$H$ is a $(-2)$-curve $\Sigma_{0} $ in the class
$h-e_{-1}-e_{0}-e_{1}$. The $\cnums ^{*}$ action extends to this blow-up
acting with two fixed points on each of the curves $E_{-1}$, $E_{0}$, and
$E_{1}$, namely the intersection with $\Sigma_{0} $ and one other. Blow-up the
fixed points on $E_{-1}$ and $E_{1}$ that are not the ones on $\Sigma_{0} $ to
obtain two new exceptional curves $E_{-2}$ and $E_{2}$ in the classes
$e_{-2}$ and $e_{2}$. Let $\Sigma _{-1}$ and $\Sigma _{1}$ be the proper
transforms of $E_{-1}$ and $E_{1}$ and note that they are $(-2)$-spheres in
the classes $e_{-1}-e_{-2}$ and $e_{1}-e_{2}$ respectively. The $\cnums
^{*}$ action extends to this blow-up and we can repeat the procedure $a-1$
additional  
times to obtain $\til{P}$. $\til{P}$ contains $2a+1$ $(-2)$-spheres, namely
$\Sigma _{-a},\ldots,\Sigma _{a}$ which represent the classes
$$[\Sigma _{n}] =\begin{cases}
e_{n}-e_{n+1}&\text{ if $0<n\leq a$,}\\
h-e_{0}-e_{-1}-e_{1}&\text{ if $a=0$,}\\
e_{n}-e_{n-1}&\text{ if $-a\leq n<0$.}
\end{cases}$$
We rename the $(-1)$-spheres $E_{0}$, $E_{a+1}$, and $E_{-a-1}$ by $\Sigma
_{*}$, $\Sigma _{a+1}$, and $\Sigma _{-a-1}$ and it is a straight forward
computation  to check that
the classes $[\Sigma _{*}],[\Sigma _{-a-1}],\ldots,[\Sigma _{a+1}]$ form an
integral basis for $H_{2}(\til{P};\znums )$.

The configuration $\Sigma _{*}+\sum_{n=-a}^{a}\Sigma _{n}$ is (as our
notation suggests) biholomorphic to $\Sigma (a)$.
Furthermore, $T\til{P}|_{\Sigma (a)}$ is isomorphic to the bundle $T$
defining $[T^{1}\to T^{2}]_{s(a)}$. This will allow us to realize our
obstruction 
problem as an ordinary Gromov-Witten invariant:
\begin{lem}\label{lem: obstr problem is same as certain GW inv in blown up P2}
$[\M _{s(a)}]^{vir}$ is the same
as the (ordinary) genus 0 Gro\-mov-Wit\-ten invariant of $\til{P} $ in the
class  
$$
[\Sigma _{*}]+\sum_{n=-a}^{a}s_{n}[\Sigma _{n}] .
$$
\end{lem}
{\sc Proof:} This follows immediately if we can show that {\em all} the
rational curves in the above homology class lie in the 
configuration $\Sigma(a) $.

Note that the curves $\Sigma _{*},\Sigma _{-a-1},\ldots,\Sigma _{a+1}$ are
preserved by the $\cnums ^{*}$ action and the only other curves preserved
are the proper transforms of lines through the fixed point $p$. We call
these  additional lines $\Sigma ^{+}$, $\Sigma ^{-}$, $\Sigma
^{0}$, and $\Sigma ^{t}$ which are the proper transforms of the lines
$\overline{pp_{+}}$, $\overline{pp_{-}}$, $\overline{pp_{0}}$, and
$\overline{pp_{t}}$ where $p_{t} $ is any point on $H$ that is
not $p_{+}$, $p_{-}$, or $p_{0}$. See the figure:
\vskip 1.25in
\setlength{\unitlength}{0.00083333in}
\begingroup\makeatletter\ifx\SetFigFont\undefined
% extract first six characters in \fmtname
\def\x#1#2#3#4#5#6#7\relax{\def\x{#1#2#3#4#5#6}}%
\expandafter\x\fmtname xxxxxx\relax \def\y{splain}%
\ifx\x\y   % LaTeX or SliTeX?
\gdef\SetFigFont#1#2#3{%
  \ifnum #1<17\tiny\else \ifnum #1<20\small\else
  \ifnum #1<24\normalsize\else \ifnum #1<29\large\else
  \ifnum #1<34\Large\else \ifnum #1<41\LARGE\else
     \huge\fi\fi\fi\fi\fi\fi
  \csname #3\endcsname}%
\else
\gdef\SetFigFont#1#2#3{\begingroup
  \count@#1\relax \ifnum 25<\count@\count@25\fi
  \def\x{\endgroup\@setsize\SetFigFont{#2pt}}%
  \expandafter\x
    \csname \romannumeral\the\count@ pt\expandafter\endcsname
    \csname @\romannumeral\the\count@ pt\endcsname
  \csname #3\endcsname}%
\fi
\fi\endgroup
{\renewcommand{\dashlinestretch}{30}
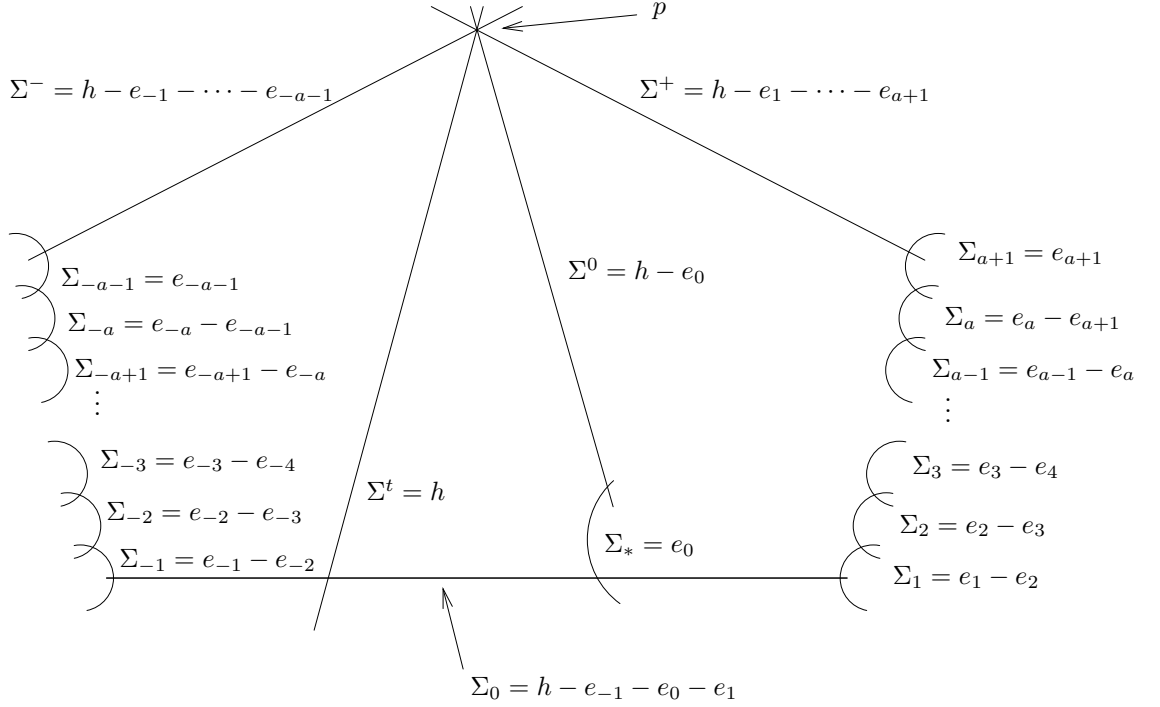
\begin{picture}(6733,4347)(0,-10)
\put(49.077,2698.116){\arc{412.042}{4.5117}{7.6590}}
\put(88.536,2373.307){\arc{413.025}{4.5198}{7.6518}}
\put(170.500,2049.500){\arc{413.021}{4.5150}{7.6566}}
\put(291.500,1401.500){\arc{413.021}{4.5150}{7.6566}}
\put(372.808,1076.961){\arc{414.198}{4.5140}{7.6537}}
\put(453.035,752.810){\arc{413.812}{4.5226}{7.6547}}
\put(5802.849,2698.228){\arc{411.851}{1.7606}{4.9136}}
\put(5763.578,2373.383){\arc{413.219}{1.7734}{4.9093}}
\put(5681.500,2049.500){\arc{413.021}{1.7682}{4.9098}}
\put(5560.500,1401.500){\arc{413.021}{1.7682}{4.9098}}
\put(5479.115,1077.077){\arc{414.002}{1.7658}{4.9112}}
\put(5399.077,752.885){\arc{414.003}{1.7706}{4.9065}}
\put(4108.585,991.708){\arc{991.136}{2.1999}{3.9794}}
\path(3210,4320)(130,2739)
\path(616,753)(5236,753)
\path(2642,4320)(5722,2739)
\path(5236,753)(616,753)
\path(2886,4320)(3777,1199)
\path(2966,4320)(1912,428)
%\path(1791,1077)(1954,793)
%\path(1868.247,882.143)(1954.000,793.000)(1920.285,912.010)
\path(2716,686)(2841,186)
\path(2774.209,576.859)(2716.000,686.000)(2716.000,562.307)
\path(3091,4186)(3925,4270)
\path(3207.390,4227.874)(3091.000,4186.000)(3213.402,4168.177)
\put(2237,1239){\makebox(0,0)[lb]{\smash{{{\SetFigFont{10}{8.4}{rm}$\Sigma ^{t}=h$}}}}}
\put(3720,915){\makebox(0,0)[lb]{\smash{{{\SetFigFont{10}{8.4}{rm}$\Sigma _{*}=e_{0}$}}}}}
\put(5925,2739){\makebox(0,0)[lb]{\smash{{{\SetFigFont{10}{8.4}{rm}$\Sigma _{a+1}=e_{a+1}$}}}}}
%\put(6127,2739){\makebox(0,0)[lb]{\smash{{{\SetFigFont{10}{8.4}{rm}$\{$a+1$\}$=e(a+1)}}}}}
\put(5844,2334){\makebox(0,0)[lb]{\smash{{{\SetFigFont{10}{8.4}{rm}$\Sigma _{a}=e_{a}-e_{a+1}$}}}}}
\put(5764,2009){\makebox(0,0)[lb]{\smash{{{\SetFigFont{10}{8.4}{rm}$\Sigma _{a-1}=e_{a-1}-e_{a}$}}}}}
\put(5844,1725){\makebox(0,0)[lb]{\smash{{{\SetFigFont{10}{8.4}{rm}$\vdots$}}}}}
\put(5641,1401){\makebox(0,0)[lb]{\smash{{{\SetFigFont{10}{8.4}{rm}$\Sigma _{3}=e_{3}-e_{4}$}}}}}
\put(5561,1036){\makebox(0,0)[lb]{\smash{{{\SetFigFont{10}{8.4}{rm}$\Sigma _{2}=e_{2}-e_{3}$}}}}}
\put(5520,712){\makebox(0,0)[lb]{\smash{{{\SetFigFont{10}{8.4}{rm}$\Sigma _{1}=e_{1}-e_{2}$}}}}}
\put(697,820){\makebox(0,0)[lb]{\smash{{{\SetFigFont{10}{8.4}{rm}$\Sigma _{-1}=e_{-1}-e_{-2}$}}}}}
\put(576,1442){\makebox(0,0)[lb]{\smash{{{\SetFigFont{10}{8.4}{rm}$\Sigma _{-3}=e_{-3}-e_{-4}$}}}}}
\put(616,1117){\makebox(0,0)[lb]{\smash{{{\SetFigFont{10}{8.4}{rm}$\Sigma _{-2}=e_{-2}-e_{-3}$}}}}}
\put(414,2009){\makebox(0,0)[lb]{\smash{{{\SetFigFont{10}{8.4}{rm}$\Sigma _{-a+1}=e_{-a+1}-e_{-a}$}}}}}
\put(373,2293){\makebox(0,0)[lb]{\smash{{{\SetFigFont{10}{8.4}{rm}$\Sigma _{-a}=e_{-a}-e_{-a-1}$}}}}}
\put(332,2576){\makebox(0,0)[lb]{\smash{{{\SetFigFont{10}{8.4}{rm}$\Sigma _{-a-1}=e_{-a-1}$}}}}}
\put(535,1766){\makebox(0,0)[lb]{\smash{{{\SetFigFont{10}{8.4}{rm}$\vdots$}}}}}
%\put(1751,450){\makebox(0,0)[lb]{\smash{{{\SetFigFont{10}{8.4}{rm}$t$}}}}}
\put(10,3752){\makebox(0,0)[lb]{\smash{{{\SetFigFont{10}{8.4}{rm}$\Sigma ^{-}=h-e_{-1}-\cdots -e_{-a-1}$}}}}}
\put(3494,2617){\makebox(0,0)[lb]{\smash{{{\SetFigFont{10}{8.4}{rm}$\Sigma ^{0}=h-e_{0}$}}}}}
\put(3940,3752){\makebox(0,0)[lb]{\smash{{{\SetFigFont{10}{8.4}{rm}$\Sigma ^{+}=h-e_{1}-\cdots -e_{a+1}$}}}}}
\put(4023,4272){\makebox(0,0)[lb]{\smash{{{\SetFigFont{10}{8.4}{rm}$p$}}}}}
\put(2885,19){\makebox(0,0)[lb]{\smash{{{\SetFigFont{10}{8.4}{rm}$\Sigma _{0}=h-e_{-1}-e_{0}-e_{1}$}}}}}
\end{picture}
}
%\vskip .5in

We express the classes of $\Sigma ^{+}$, $\Sigma
^{-}$, $\Sigma ^{0}$ and $\Sigma ^{t}$ in the two homology bases
$\{h,e_{-a-1},\ldots,e_{a+1} \}$ and $\{[\Sigma _{*}],[\Sigma
_{-a-1}],\ldots,[\Sigma _{a+1}]  \}$ as follows:
\begin{eqnarray*} 
\left[ \Sigma  ^{t}\right] &=&h=\left[ \Sigma
_{*}\right] +\sum_{n=-a-1}^{a+1}\left[ \Sigma _{n}\right] \\ 
\left[  \Sigma ^{0} \right] &=&h-e_{0}=\sum_{n=-a-1}^{a+1}\left[ \Sigma
_{n}\right] \\ 
\left[ \Sigma ^{+} \right] &=&h-e_{1}-\cdots -e_{a+1}\\
&=&\left[ \Sigma _{*}\right] +\sum_{n=-a-1}^{a+1}\left[ \Sigma
_{n}\right] -\sum_{n=1}^{a+1}n\left[ \Sigma _{n}\right]  \\
\left[ \Sigma ^{-}\right] &=&h-e_{-1}-\cdots -e_{-a-1}\\
&=&\left[ \Sigma _{*}\right] +\sum_{n=-a-1}^{a+1}\left[ \Sigma
_{n}\right] -\sum_{n=1}^{a+1}n\left[ \Sigma _{-n}\right] . 
\end{eqnarray*}

Since $\cnums ^{*}$ acts on $\til{P}$ we get an $\cnums ^{*}$-action on the
moduli space 
$\M ^{\til{P}}_{\{s_{n} \}}$ of genus 0 stable maps to $\til{P}$ in the
class $[\Sigma _{*}]+\sum_{n=-a}^{a}s_{n}[\Sigma _{n}]$. We first show that
the maps in the fixed point set of $\M ^{\til{P}}_{\{s_{n} \}}$ must have
image $\Sigma _{*}+\sum_{n=-a}^{a}s_{n}\Sigma _{n}$. This is essentially
for homological reasons: the image of a map in the fixed point set of $\M
^{\til{P}}_{\{s_{n} \}}$ must be of the form 
$$c_{*}\Sigma _{*}+\sum_{n=-a-1}^{a+1}c_{n}\Sigma _{n}+\sum_{t} c^{t}\Sigma
^{t}+c^{+}\Sigma ^{+}+c^{-}\Sigma ^{-} +c^{0}\Sigma ^{0}$$  
for non-negative coefficients given by the $c$'s. Since $[\Sigma
_{*}],[\Sigma _{-a-1}],\ldots,[\Sigma _{a+1}]$ form a basis we have
\begin{eqnarray*}
c_{*}+\sum_{t}c^{t}+c^{+}+c^{-}&=&1,\\
c_{n}+c^{0}+\sum_{t}c^{t}+c^{+}+c^{-}-|n|c^{\operatorname{sign}(n)}&=
&\begin{cases} 
s_{n}&|n|\leq a,\\
0&|n|=a+1.
\end{cases}
\end{eqnarray*}

The first equation implies that exactly one of $c_{*}$, $c^{t}$, $c^{+}$,
or $c^{-}$ is 1 (for some $t$) and the rest are 0. Suppose that $c^{+}=1$;
then $c^{-}=0$ 
and letting $n=-a-1$ in the second equation leads to a contradiction and so
we have $c^{+}=0$. A similiar argument shows $c^{-}=0$ and then summing the
second equation over $n$ leads to 
$$\left(\sum_{n=-a-1}^{a+1}c_{n} \right)+(2a+3)(c^{0}+\sum_{t}c^{t})=a $$
which implies that $c^{0}=c^{t}=0$. Thus $c_{*}=1$ and $c_{n}=s_{n}$.

Finally, suppose $f\in \M ^{\til{P}}_{\{s_{n} \}}$ is not a fixed point of
the $\cnums ^{*}$ action. Then the limit of the action of $\lambda \in
\cnums ^{*}$ on $f$ as $\lambda \to 0$ must be fixed an hence has image
$\Sigma _{*}+\sum_{n=-a}^{a}s_{n}\Sigma _{n}$. But then the limit of the
action as $\lambda \to \infty $ must also be fixed and its image must
contain the point $p$ which 
is a contradiction. Hence every $f\in \M ^{\til{P}}_{\{s_{n} \}}$ is fixed
by $\cnums ^{*}$ and so has image
$\Sigma _{*}+\sum_{n=-a}^{a}s_{n}\Sigma _{n}$. \qed

Now $\til{P}$ is deformation equivalent to the blow-up of $\P ^{2}$ at $2a+3$
generic points and so the invariant for the  class $[\Sigma
_{*}]+\sum_{n=-a}^{a}s_{n}[\Sigma _{n}]$  can be computed
using elementary properties of the invariants for blow-ups of $\P ^{2}$. 

We follow the notation of \cite{Go-Pa} and recall some of the properties of
the invariant. We write $N{(d;\alpha
_{1},\ldots)}$ for the genus 0 Gro\-mov-Wit\-ten invariant in the class
$dh-\sum_{i}\alpha _{i}e_{i}$ . Here we are not being very picky about the
indexing set for the exceptional classes since the invariant is the same
under reordering. In the notation $N{(d;\alpha
_{1},\ldots)}$ it is implicit that if the moduli space of genus 0 maps in
the class $(d;\alpha _{1},\ldots)$ is positive dimensional then we impose
the proper number of point constraints and if the dimension is negative the
invariant is zero. Also, we drop any $\alpha =0$ terms from the notation so
that $N(d;\alpha _{1},\ldots,a_{l},0,\ldots,0)=N(d;\alpha
_{1},\ldots,\alpha _{l})$.  The invariants satisfy the following properties:
\begin{enumerate}
\item $N{(d;\alpha _{1},\ldots)}=0$ if any $\alpha <0$ unless $d=0$, $\alpha
_{i}=0$ for all $i$ except $i_{0}$ and $\alpha _{i_{0}}=-1$. In the latter
case the invariant is 1.
\item $N{(d;\alpha _{1},\ldots,\alpha _{l},1)}=N{(d;\alpha
_{1},\ldots,\alpha _{l} )}$.
\item $N{(d;\alpha _{1},\ldots,\alpha _{l})}=N{(d;\alpha _{\sigma
(1)},\ldots,\alpha _{\sigma (l)})}$ for any permutation $\sigma $.
\item $N{(d;\alpha _{1},\ldots,\alpha _{l})}$ is invariant under the
Cremona transformation which takes the class
 $$(d;\alpha _{1},\alpha _{2},\alpha
_{3},\ldots)$$
 to the class 
$$(2d-\alpha _{1}-\alpha _{2}-\alpha _{3};d-\alpha
_{2}-\alpha _{3},d-\alpha _{1}-\alpha _{3},d-\alpha _{1}-\alpha
_{2},\ldots).$$
\item $N{(1)}=1$.
\end{enumerate}

Ordering the exceptional classes in $\til{P}$ by
$e_{0},e_{1},e_{-1},e_{2},e_{-2},\ldots$ and rewriting the class $\Sigma
_{*}+\sum 
s_{n}\Sigma _{n}$ in this basis, we can express can express the
contribution of $\M _{s(a)}$ as
$$
[\M_{s(a)}]^{vir}=
N{(s_{0};s_{0}-1,s_{0}-s_{1},s_{0}-s_{-1},s_{1}-s_{2},s_{-1}-s_{-2},
\ldots,s_{-a+1}-s_{-a})}. 
$$

We call an admissible sequence $\{s_{n} \}$ {\em 1-admissible} if
$s_{\pm n\pm 1}$ is either $s_{\pm n}$ or $s_{\pm n}-1$ for all $n$.

\begin{lem}\label{lem: contribution is 1 for 1-admissable, 0 otherwise}
$[\M _{s(a)}]^{vir}=1$ if $s(a)$ is a 1-admissible sequence and $[\M
_{s(a)}]^{vir}=0$ otherwise.
\end{lem}
{\sc Proof:}
Suppose that  $[\M_{s(a)}]^{vir}\neq 0$. Since $s_{0}>0$, all the other
terms in 
$(s_{0};s_{0}-1,s_{0}-s_{1},\ldots)$ must be non-negative by property 1. Thus
$s_{\pm n\pm 1}\geq s_{\pm n}$ for all $n$. Now by permuting and performing the
Cremona transformation, we get
\begin{eqnarray*}
[\M_{s(a)}]^{vir}&=&N{(s_{0};s_{0}-1,s_{0}-s_{1},s_{\pm n}-s_{\pm n\pm
1},\ldots)}\\  
&=&N(1+s_{1}+s_{\pm n\pm 1}-s_{\pm n};s_{1}-s_{\pm n}+s_{\pm n\pm 1}, \\
& &\quad 
1+s_{\pm n\pm 1}-s_{\pm n},s_{1}+1-s_{0},\ldots).
\end{eqnarray*}
Now since $s_{\pm n}\leq s_{1}$, we have $1+s_{1}+s_{\pm n\pm 1}-s_{\pm
n}>0$ and so 
$1+s_{\pm n\pm 1}-s_{\pm n}\geq 0$ which combined with $s_{\pm n\pm 1}\leq
s_{\pm n}$ 
yields
$$s_{\pm n\pm 1}\leq s_{\pm n}\leq s_{\pm n\pm 1} $$ 
and so $s$ is 1-admissible.

Suppose then that $s$ is 1-admissible. Then except for the first two terms,
the class
$(s_{0};s_{0}-1,s_{0}-s_{1},\ldots)$ consists of $0$'s and $1$'s. Thus
$[\M_{s(a)}]^{vir}=N{(s_{0};s_{0}-1)}$. Finally, since
$N{(s_{0};s_{0}-1)}=N{(s_{0};s_{0}-1,1,1)}$ we can apply Cremona to get 
$$N{(s_{0};s_{0}-1)}=N{(s_{0}-1;s_{0}-2)} $$
and so by induction
$$[\M_{s(a)}]^{vir}=N{(s_{0};s_{0}-1)}=N{(1)}=1 $$
and the lemma is proved.\qed 

\begin{lem}\label{lem: no. of 1-admissable seqs is p(a)}
The number of 1-admissible sequences $s$ with $|s|=a$ is the number of
partitions of $a$, $p(a)$.
\end{lem}
{\sc Proof:}\footnote{We are grateful to D. Maclagan and S. Schleimer for
help with this and other combinatorial difficulties.} 
The number of partitions, $p(a)$, is given by the number of Young diagrams
of size $a$. There is a bijective correspondence between 1-admissible
sequences and Young diagrams. Given a Young diagram define an 1-admissible
sequence $\{s_{n} \}$ by setting $s_{0}$ equal to the number of blocks on
the diagonal, $s_{1}$ equal to the number of blocks on the first lower
diagonal, $s_{2}$ equal to the number of blocks on the second lower
diagonal, and so on, doing the same for $s_{-1}$, $s_{-2}$,\dots with the
upper diagonals. It is easily seen that this defines a bijection.\qed

Summarizing the results of this section we have:
\begin{thm}\label{thm: contribution of Mab}
Since every component of $\M _{\mathbf{a},\mathbf{b}}$ contributes either 0
or 1 to $N_{g}(n)$, the overall contribution of  $\M
_{\mathbf{a},\mathbf{b}}$ is the sum 
over all the connected components whose contribution is 1. It is thus the
sum of all choices of data $\Lambda (\mathbf{b})$ and $\lambda
(\mathbf{b})$, and those choices of $s(\mathbf{a})$ that are 1-admissible. 
For each $a_{j}\in \mathbf{a}$ we have $p(a_{j})$ choices of a 1-admissible
sequence $s(a_{j})$ and for each $b_{i}\in \mathbf{b}$ we have $b_{j}\sigma
(b_{i})$ choices for the data $\Lambda (b_{i})$ and $\lambda (b_{i})$. Thus
the total contribution is:
$$[\M
_{\mathbf{a},
\mathbf{b}}]^{vir}=\prod_{j=1}^{24}p(a_{j})\prod_{i=1}^{g}b_{i}\sigma 
(b_{i}). $$  
\end{thm}

\section{Counting curves on the rational elliptic surfaces}\label{sec:
counting on E1} 

Let $Y$ be the blow up of
$\mathbf{P}^{2}$ at nine 
distinct points. In this section we apply our degeneration method and our
local calculations to compute a certain set of Gro\-mov-Wit\-ten invariants of
$Y$. We compute the genus $g$ invariants for all classes such that the
invariants require exactly
$g$ constraints.  There is a canonical symplectic form $\omega $ (unique up
to deformation equivalence) on $Y$ determined by the blow up of the
Fubini-Study form on $\mathbf{P}^{2}.$    
                           
If we arrange these nine blow up points lying on a pencil of cubic elliptic
curves in $\mathbf{P}^{2},$ then $Y$ has the structure of an elliptic surface
with fiber class $F$ representing these elliptic curves in $H_{2}\left(
Y,\mathbf{Z}\right) $ and the nine exceptional curves $e_{1},e_{2},...,e_{9}$
are all sections of this elliptic fibration. If $h$ represents the homology
class of the strict transform of the hyperplane in $\mathbf{P}^{2},$ then we
have $F=3h-e_{1}-...-e_{9}$. In fact $H^{2}\left( Y,\mathbf{Z}\right) $ is
generated by $e_{1},...,e_{9}$ and $h$. We abbreviate the class
$dh-a_{1}e_{1}-\cdots -a_{9}e_{9}$ by $(d;a_{1},\ldots,a_{9})$. 

Now we pick any of these sections, $e_{9}$ say, and consider the class
$C_{n}= e_{9}+( g+n) F=(3(n+g);g+n,\ldots,g+n,g+n-1)$. It is easy to check
that the complete linear system $|C_{n}|$ has dimension $g+n.$ We write
$N_{g}^{Y}(C ) $ the Gro\-mov-Wit\-ten invariant for $\left( Y,\omega \right) $
which counts the number of curves of geometric genus $g$
representing the homology class $C$ and passing
through $g$ points. \ie we define 
$$N_{g}^{Y}(C)=\Psi
^{Y}_{(C,g,g)}(PD(x_{1}),\ldots,PD(x_{g});PD(\overline{\M }_{g,g})) . $$

We show that the numbers $N_{g}^{Y}(C_{n})$ contain all the genus $g$
Gro\-mov-Wit\-ten invariants that are constrained to exactly $g$ points. This
was observed by G\"ottsche who explained the following argument to us:

For $N_{g}^{Y}(C)$ to be well defined (see Equation \ref{eqn:
dimension formula}) we need $4g=2c_{1}(Y)\cdot
C+2g-2(1-g)$, \ie $F\cdot C=1$. Now the Gro\-mov-Wit\-ten invariants do not
change when $C\mapsto C'$ is induced by a  permutation of the exceptional
classes $e_{i}$ or a Cremona transform (see \cite{Go-Pa}). Recall that the
Cremona 
transform takes a class 
$(d;a_{1},\ldots,a_{9})$ to the class 
$$
(2d-a_{1}-a_{2}-a_{3};d-a_{2}-a_{3},
d-a_{1}-a_{3},d-a_{1}-a_{2},a_{4},\ldots,a_{9}). 
$$
\begin{lem}\label{lem: lothar's lemma}
Let $C\in H_{2}(Y;\znums )$ be a class so that the moduli space of genus
$g$ maps has formal dimension $g$. Then the
class $C$ can be transformed by a
sequence of Cremona transforms and permutations of the $e_{i}$'s to a class
of the form $e_{9}+(g+n)F=C_{n}$. 
\end{lem}
{\em Proof:} By permuting the $E_{i}$'s we may assume that
$a_{1}\geq a_{2}\geq \cdots \geq a_{9}$. Then the condition $F\cdot C=1$ is
equivalent to $3d-1=\sum_{i}a_{i}$ so that $a_{1}+a_{2}+a_{3}\geq d$ with
equality if and only if $C=(3i,i,i,i,i,i,i,i,i,i-1)=e_{9}+iF$ for some
$i=n+g$. If the equality is strict then we can apply a Cremona transform to
obtain $C'=(e,b_{1},\ldots,b_{9})$ with $e<d$. The result follows by
descending induction on $d$.\qed 

The methods of Sections \ref{sec: computation of Ngn} and \ref{sec:
analysis and local computations} apply to these invariants (see remark
\ref{rem: obstr argumnet applies to E(1) also}). Note that the elliptic
fibration of $Y$ has (generically) 12 nodal fibers rather than 24. We
get the same formula as in the $K3$ case with the 24 replace by 12.
\begin{theorem}\label{thm: E1 calculation}
For any $g\geq 0,$ we have
\begin{eqnarray*}
\sum_{n=0}^{\infty }N_{g}^{Y}(C_{n}) q^{n}
&=&\left(\sum_{b=1}^{\infty }b\sigma (b)q^{b-1} \right)^{g}\prod_{m=1}^{\infty
}\left( 1-q^{m}\right) ^{-12} \\ 
&=&\left(\frac{d}{dq}G_{2}(q) \right)^{g}\left(\frac{q}{\Delta (q)} \right)^{1/2}
\end{eqnarray*}  
\end{theorem}

When the genus $g$ equals zero, these numbers are computed by G\"{o}ttsche
and Pandharipande \cite{Go-Pa}. In fact, they obtain
all genus zero Gro\-mov-Wit\-ten invariants for $\mathbf{P}^{2}$ blown up at
arbitrary number of points in terms of two rather complicated recursive
formulas. Theorem \ref{thm: E1 calculation} can be verified term by term
for $g=0$ 
using the recurrence relations, although the computer calculation becomes
extremely lengthy quickly. We know of no way of obtaining the genus 0
closed form 
of Theorem \ref{thm: E1 calculation} directly from the recurrence relations.

%\bibliography{mainbiblio}
%\bibliographystyle{plain}

\end{document}